\pdfoutput=1
\documentclass[toc]{mathprint}
\usepackage{rutarmacros}
\usepackage{macros}

\title[Attainable dimensions]{Attainable forms of intermediate dimensions}

\author[Banaji]
  {Amlan Banaji}
  {Mathematical Institute, University of St Andrews, St Andrews KY16 9SS, Scotland}
  {afb8@st-andrews.ac.uk}

\author[Rutar]
  {Alex Rutar}
  {Mathematical Institute, University of St Andrews, St Andrews KY16 9SS, Scotland}
  {alex@rutar.org}

\journalcite{zbl:1509.28005}

\begin{document}
\begin{abstract}
    The intermediate dimensions are a family of dimensions which interpolate between the Hausdorff and box dimensions of sets.
    We prove a necessary and sufficient condition for a given function $h(\theta)$ to be realized as the intermediate dimensions of a bounded subset of $\mathbb{R}^d$.
    This condition is a straightforward constraint on the Dini derivatives of $h(\theta)$, which we prove is sharp using a homogeneous Moran set construction.
\end{abstract}

\section{Introduction}
The Hausdorff and box dimensions of a set $F\subset\R^d$ are two widely studied notions of dimension.
We denote these dimensions by $\dimH F$ and $\dimB F$ respectively.
The box dimension is a coarse measurement of dimension, in the sense that it only takes into account the size of the set at a fixed scale $\delta$, as $\delta$ goes to zero.
On the other hand, the Hausdorff dimension takes into account all small scales simultaneously.
When the box and Hausdorff dimensions are equal, they indicate that the set has a large amount of spatial regularity.
For example, in \cite{shm2019b}, Shmerkin established Falconer's distance problem for subsets of the plane with equal Hausdorff and box dimensions, whereas this conjecture is wide open in general.

However, the box and Hausdorff dimensions of sets can also differ: this is often the case for natural ``fractal'' sets, such as self-affine sets.
In order to better understand this situation, the authors of \cite{ffk2020} introduced the \defn{intermediate dimensions} of the set $F$.
This is a family of dimensions, which we will denote by $\dim_\theta F$ for $\theta\in[0,1]$, which satisfy $\dim_0 F=\dimH F$, $\dim_1 F=\dimB F$, and which are continuous on the interval $(0,1]$ (though not necessarily at 0).
The intermediate dimensions are defined by allowing an increasing amount of flexibility in the sizes of the covering sets as $\theta$ tends to zero.
Precise definitions are given in the next section.
We remark that intermediate dimensions are an example of dimension interpolation, of which the Assouad spectrum is also an example.
We refer the reader to \cite{fra2021a} for a recent survey of this topic.

In a very heuristic sense, the intermediate dimensions for $\theta>0$ behave more like box dimensions than Hausdorff dimensions.
Assuming the intermediate dimensions are continuous at zero, one would hope to obtain information about the Hausdorff dimension of a set $F$ in terms of the intermediate dimensions, which may \emph{a priori} be easier to bound.
For example, in \cite{bff2021}, as an application of general results on intermediate dimensions of projections, the authors obtain novel bounds on the box dimensions of projections of sets with intermediate dimensions continuous at zero.
Continuity at zero also provides information about the box dimensions of images of sets under fractional Brownian motion \cite{bur2021}.

Moreover, intermediate dimensions can distinguish bi-Lipschitz equivalence even when other notions of dimension give no information \cite{bk2021}.
In fact, quantitative information about Hölder exponents can be obtained for any bi-Hölder map between sets with distinct intermediate dimensions \cite[Theorem~3.1]{bur2021}.
We refer the reader to \cite{fal2021} for a survey of some recent results concerning intermediate dimensions.

In this paper, we study the forms of the intermediate dimensions $\dim_\theta F$ for $\theta\in[0,1]$.
In general, the intermediate dimensions are known to satisfy certain regularity constraints (see \cite[Prop~2.1]{ffk2020} and \cite[Proposition~2.1]{fal2021}).
On the other hand, intermediate dimensions have been computed for some specific families of sets.
For example, Bedford--McMullen carpets are a recent example of a natural family of sets for which the intermediate dimensions exhibit interesting properties \cite{bk2021}.
Other sets which have been studied in the literature include infinitely generated self-conformal sets \cite{bf2023} and elliptical polynomial spirals \cite{bff2021}.
However, in general, no progress has been made on determining sharpness of the general constraints on the intermediate dimensions.

Our main result is to obtain a full characterization of possible intermediate dimensions for subsets of $\R^d$.

Given a function $f\colon\R\to\R$, we denote the (upper right) Dini derivative (see \cref{e:Dini-intro}) of $f$ at $x$ by $\diniu{+}f(x)$.
We then have the following result.
\begin{itheorem}\label{it:intro-main-res}
    Let $h\colon[0,1]\to[0,d]$ be any function.
    Then there exists a non-empty bounded set $F\subset\R^d$ with $\dim_\theta F=h(\theta)$ if and only if $h$ is non-decreasing, is continuous on $(0,1]$, and satisfies
    \begin{equation}\label{e:intro-gen-bound}
        \diniu{+}h(\theta)\leq\frac{h(\theta)(d-h(\theta))}{d\theta}
    \end{equation}
    for all $\theta\in(0,1)$.
\end{itheorem}
We see that the intermediate dimensions can have highly varied behaviour; such behaviour has not been seen in any prior examples.
In particular, without stronger assumptions on the set $F$, very little can be said about the possible forms of the intermediate dimensions.

For example, it follows directly from \cref{e:intro-gen-bound} that if $f$ is any non-decreasing Lipschitz function on $[0,1]$, there exists some constants $a>0$, $b\in\R$, and a set $F\subset\R$ such that $\dim_\theta F=af(\theta)+b$ for all $\theta\in[0,1]$.
In particular, the following behaviours for the intermediate dimensions are all possible:
\begin{enumerate}[nl,r]
    \item Constant on countably many disjoint closed intervals in $[0,1]$, and strictly increasing otherwise.
    \item Strictly concave, strictly convex, or linear and non-constant, on $[0,1]$.
    \item Non-differentiable at each point in a dense subset $E$ of $(0,1)$ with $\dimH E=1$ (in fact, the points of non-differentiability can form an arbitrary $G_{\delta\sigma}$ subset of $[0,1]$ with Lebesgue measure zero~\cite{zah1946}).
\end{enumerate}
This resolves all remaining questions asked in Falconer's survey \cite{fal2021}.

We observe that the bound \cref{e:intro-gen-bound} is attained at all $\theta\in(0,1)$ for the sets
\begin{equation*}
    G_{p,d}\coloneqq \bigl\{x/\norm{x}^2:x\in\{n^p:n\in\N\}^d\bigr\}
\end{equation*}
where $d\in\N$ and $p>0$.
These sets have intermediate dimensions given by $\dim_\theta G_{p,d}=d\theta/(p+\theta)$ for all $\theta\in[0,1]$ \cite[Proposition~3.8]{bf2023} (see also \cite[Prop~3.1]{ffk2020} for the case $d=1$).

Results similar to the existence result in \cref{it:intro-main-res} have also been obtained for the Assouad spectrum \cite{fhh+2019}, though a full characterization for the Assouad spectrum is not known.
Both our results, and the results in that paper, use homogeneous Moran sets as the basis of the construction.
\subsection{Notation}
We fix some $d\in\N$ and work in $\R^d$, equipped with the max norm.
Given $x\in\R^d$, we denote the $j$\textsuperscript{th} coordinate of $x$ by $x^{(j)}$.
We write $B(x,r)$ to denote the open ball with radius $r$ centred at $x$.

All sets $F$ are non-empty bounded subsets of $\R^d$.
We write $\overline{F}$ to denote the topological closure of $F$.
We also denote by $N_r(F)$ the minimal number of sets with diameter $r$ required to cover $F$.

\subsection{Statement and summary of main results}
We begin with a precise definition of the intermediate dimensions.
\begin{definition}
    Let $F\subset\R^d$.
    For $0\leq\theta\leq 1$, the \defn{upper intermediate dimensions} are given by
    \begin{align*}
        \overline{\dim}_\theta F = \inf\Bigl\{s\geq 0:&(\forall\epsilon>0)\,(\exists\delta_0>0)\,(\forall0<\delta\leq\delta_0)\,(\exists\text{ a cover }\{U_i\}_{i=1}^\infty\text{ of }F)\text{ s.t. }\\*
                                                      &\delta^{1/\theta}\leq\diam U_i\leq\delta\text{ and }\sum_{i=1}^\infty(\diam U_i)^s\leq\epsilon\Bigr\}.
    \end{align*}
    Similarly, the \defn{lower intermediate dimensions} are given by
    \begin{align*}
        \underline{\dim}_\theta F = \inf\Bigl\{s\geq 0:&(\forall\epsilon>0)\,(\forall\delta_0>0)\,(\exists<\delta\leq\delta_0\text{ and a cover }\{U_i\}_{i=1}^\infty\text{ of }F)\text{ s.t. }\\*
                                                      &\delta^{1/\theta}\leq\diam U_i\leq\delta\text{ and }\sum_{i=1}^\infty(\diam U_i)^s\leq\epsilon\Bigr\}.
    \end{align*}
\end{definition}
Clearly $\underline{\dim}_\theta F\leq\overline{\dim}_\theta F$.
If $\underline{\dim}_\theta F=\overline{\dim}_\theta F$, we denote this common value by $\dim_\theta F$.
We make the following basic observations:
\begin{itemize}
    \item $\dimH F=\underline{\dim}_0 F=\overline{\dim}_0 F$,
    \item $\underline{\dim}_\theta F$ and $\overline{\dim}_\theta F$ are monotonically increasing in $\theta$, and
    \item $\underline{\dim}_1 F=\dimlB F$ and $\overline{\dim}_1 F=\dimuB F$.
\end{itemize}
We also recall the definitions of the Assouad and lower dimensions of the set $F$.
As we will see, these dimensions influence the possible forms of the intermediate dimensions in a natural way.
\begin{definition}
    Let $F\subset\R^d$.
    The \defn{Assouad dimension} of $F$ is given by
    \begin{align*}
        \dimA F = \inf\Bigl\{\alpha:&(\exists C>0)\,(\forall 0<r<R\text{ and }x\in F),\\*
                                    &N_r(B(x,R)\cap F)\leq C\left(\frac{R}{r}\right)^\alpha\Bigr\}
    \end{align*}
    and, dually, the \defn{lower dimension} of $F$ is given by
    \begin{align*}
        \dimL F = \sup\Bigl\{\lambda:&(\exists C>0)\,(\forall 0<r<R\leq\diam F\text{ and }x\in F),\\*
                                    &N_r(B(x,R)\cap F)\geq C\left(\frac{R}{r}\right)^\lambda\Bigr\}.
    \end{align*}
\end{definition}
In general, we have $\dimL F\leq\dimH F\leq\dimlB F\leq\dimuB F\leq\dimA F$.
We refer the reader to \cite{fal1997} and \cite{fra2021} for more details on these notions of fractal dimension.

Recall that the upper Dini derivative of a function $f\colon\R\to\R$ at $x$ is given by
\begin{equation}\label{e:Dini-intro}
    \diniu{+}f(x)=\limsup_{\epsilon\to 0^+}\frac{f(x+\epsilon)-f(x)}{\epsilon}.
\end{equation}
We then define the following classes of functions.
\begin{definition}\label{d:h-class}
    Let $0\leq\lambda\leq\alpha\leq d$.
    If $\lambda<\alpha$, we denote by $\mathcal{H}(\lambda,\alpha)$ the set of functions $h\colon[0,1]\to[\lambda,\alpha]$ which satisfy the following constraints:
    \begin{enumerate}[nl,r]
        \item $h$ is non-decreasing,
        \item $h$ is continuous on $(0,1]$, and
        \item for each $\theta\in(0,1)$, we have
            \begin{equation}\label{e:h-bound-general}
                \diniu{+}h(\theta)\leq\frac{(h(\theta)-\lambda)(\alpha-h(\theta))}{(\alpha-\lambda)\theta}.
            \end{equation}
    \end{enumerate}
    Otherwise, $\lambda=\alpha$ and we let $\mathcal{H}(\lambda,\alpha)$ be the set consisting only of the constant function $h(\theta)=\alpha$.
\end{definition}
We note that in \cref{e:h-bound-general} one could take instead the lower Dini derivative and the class of functions would remain unchanged (this follows from \cref{c:dmvt}).

We now state our main result precisely.
\begin{itheorem}\label{it:general-form}
    Suppose $F\subset\R^d$ has $\dimL F=\lambda$ and $\dimA F=\alpha$.
    Then with $\underline{h}(\theta)=\underline{\dim}_\theta F$ and $\overline{h}(\theta)=\overline{\dim}_\theta F$, we have $\underline{h},\overline{h}\in\mathcal{H}(\lambda,\alpha)$, $\underline{h}\leq\overline{h}$, and $\underline{h}(0)=\overline{h}(0)$.

    Conversely, if $0\leq\lambda\leq\alpha\leq d$ and $\underline{h},\overline{h}\in\mathcal{H}(\lambda,\alpha)$ satisfy $\underline{h}\leq\overline{h}$ and $\underline{h}(0)=\overline{h}(0)$, then there exists a compact perfect set $F\subset\R^d$ such that $\dimL F=\lambda$, $\dimA F=\alpha$, $\underline{\dim}_\theta F=\underline{h}(\theta)$, and $\overline{\dim}_\theta F=\overline{h}(\theta)$ for all $\theta\in[0,1]$.
\end{itheorem}
The proof of this result is given in \cref{p:int-dim-bound} and \cref{c:upper-lower-match}.
This result gives a full characterization of all possible forms of the upper and lower intermediate dimensions of a bounded set $F\subset\R^d$.

Constraint \cref{e:h-bound-general} generalizes all previously known bounds, namely those in \cite{fal2021,ffk2020} and the first arXiv version of \cite{ban2020}.
Note that \cref{e:h-bound-general} also provides quantitative information about the Assouad and lower dimensions in terms of the intermediate dimensions.
This is in contrast to the box and Hausdorff dimensions, which provide no more information about the Assouad and lower dimensions beyond the usual order constraints.
We can also view the bound in \cref{e:h-bound-general} as $(2\theta)^{-1}$ times the harmonic mean of $h(\theta)-\lambda$ and $\alpha-h(\theta)$.
In particular, if $0\leq\lambda'\leq\lambda\leq\alpha\leq\alpha'\leq d$, then $\mathcal{H}(\lambda',\alpha')\supseteq\mathcal{H}(\lambda,\alpha)$.
Of course, by taking $\underline{h}=\overline{h}$ we can also ensure that the intermediate dimensions exist.
Therefore, \cref{it:intro-main-res} follows from \cref{it:general-form}.

The proof of the bound \cref{e:h-bound-general} is given in \cref{p:int-dim-bound}.
The strategy is essentially as follows.
Given $\epsilon>0$, we want to convert an optimal cover $\mathcal{U}$ for some set of scales $[\delta^{1/\theta},\delta]$ into a cover for a smaller set of scales $[{\delta'}^{1/(\theta+\epsilon)},\delta']\subset[\delta^{1/\theta},\delta]$.
If $U\in\mathcal{U}$ has large diameter, then we can replace it using sets with smaller diameter using the Assouad dimension; and we can cover the $U\in\mathcal{U}$ with small diameter with fewer sets of larger diameter using the lower dimension.
Then $\delta'$ is chosen carefully to optimize this process.

In order to establish a converse to this general bound, our main strategy is to construct sets which we call \defn{homogeneous Moran sets}.
These sets are analogous to the $2^d$-corner Cantor sets in $\R^d$, except we only require the subdivision ratios to be equal within each stage in the construction, and not necessarily between stages.
The following nice property was essentially observed in \cite{chm1997}: the optimal covers for a homogeneous Moran set can be taken to consist of sets with equal diameter.
This result is given in \cref{l:flat-covers}.
A direct application of this result is a convenient formula for the upper intermediate dimensions of these sets, given in \cref{p:int-dim}.

Using this formula, in \cref{l:exact-construction} and \cref{l:Moran-formula}, we present a general strategy to construct homogeneous Moran sets with upper intermediate dimensions given by an infimum over a ``sliding window'' of a function $g$ satisfying certain derivative constraints.
Then for any $h(\theta)$ satisfying the general bounds, in \cref{t:upper-h-form} we construct a function satisfying the derivative constraints so that the corresponding Moran set has upper intermediate dimensions given by the prescribed formula.
This establishes \cref{it:general-form} for the upper intermediate dimensions.

Finally, in \cref{t:gen-h-form}, we construct an inhomogeneous Moran set which, at a fixed scale, looks like a finite union of homogeneous Moran sets each with prescribed upper intermediate dimension $h(\theta)$.
This process is done in such a way to ensure that the intermediate dimensions exist.
Then, taking a disjoint union of this set with the set provided in \cref{t:upper-h-form}, we complete the proof of \cref{it:general-form}.
The details are provided in \cref{c:upper-lower-match}.

Heuristically, the covering strategy for \cref{p:int-dim-bound} will be sharp when the relative covering numbers in the Assouad and lower dimensions are realized uniformly on the entire set for a fixed scale.
In some sense, this motivates the choice of homogeneous Moran sets, which have the maximum possible uniformity at a fixed scale.
The key observation is that inhomogeneity between scales is sufficient to obtain all possible forms of the intermediate dimensions.

\section{Intermediate dimensions and general bounds}
\subsection{Some elementary results on Dini derivatives}
We begin with some standard results on Dini derivatives, which will be useful later in the paper.
We refer the reader to \cite{bru1994} for more details.
\begin{definition}
    Let $g\colon\R\to\R$ be a function.
    Then the \defn{upper right Dini derivative} is given by
    \begin{equation*}
        \diniu{+}g(x)=\limsup_{\epsilon\to 0^+}\frac{g(x+\epsilon)-g(x)}{\epsilon}.
    \end{equation*}
    The lower right Dini derivative is denoted $\dinil{+}g$, and the left Dini derivatives are analogously denoted $\diniu{-}g$ and $\dinil{-}g$.
\end{definition}
We first make the following observation.
\begin{lemma}\label{l:c1-bound}
    Let $f$ and $g$ be continuous functions on $[a,b]$ with $\dinil{+}g\leq \diniu{+}f$ and $g(a)=f(a)$.
    Then $g\leq f$.
\end{lemma}
\begin{proof}
    Observe that $\diniu{+}(f-g)=\diniu{+}f-\dinil{+}g\geq 0$ so by \cite[Corollary~11.4.2]{bru1994}, $f-g$ is non-decreasing.
    But $(f-g)(a)=0$ so $g\leq f$.
\end{proof}
As an application, we obtain the following analogue of the mean value theorem.
\begin{corollary}\label{c:dmvt}
    Let $g$ be a continuous function on $[a,b]$ and set $s=\frac{g(b)-g(a)}{b-a}$.
    Then for any $\phi\in\{\diniu{+}g,\dinil{+}g,\diniu{-}g,\dinil{-}g\}$,
    \begin{enumerate}[nl,r]
        \item there exists $x\in[a,b]$ such that $\phi(x)\leq s$, and
        \item there exists $x\in[a,b]$ such that $\phi(x)\geq s$.
    \end{enumerate}
\end{corollary}
\begin{proof}
    We prove that there is some $x$ such that $\dinil{+}g(x)\geq s$; the other cases are similar.
    Without loss of generality, there is some $x_0\in(a,b)$ such that $g(x_0)> g(a)+s(x_0-a)$.
    Suppose for contradiction $\dinil{+}g(x)\leq s$ for all $x\in[a,x_0]$.
    By \cref{l:c1-bound}, $g(x)\leq s(x-a)+g(a)$ for all $x\in[a,x_0]$, contradicting the choice of $x_0$.
\end{proof}
We now have the following elementary result.
\begin{lemma}\label{l:g-check}
    Let $0\leq\lambda\leq\alpha\leq d$, let $g\colon\R\to(\lambda,\alpha)$ be continuous, and let $U\subset\R$ be an open set.
    Then the following are equivalent:
    \begin{enumerate}[nl,r]
        \item $\diniu{+}g(x)\in[\lambda-g(x),\alpha-g(x)]$ for all $x\in U$.
        \item $\dinil{+}g(x)\in[\lambda-g(x),\alpha-g(x)]$ for all $x\in U$.
        \item $\diniu{-}g(x)\in[\lambda-g(x),\alpha-g(x)]$ for all $x\in U$.
        \item $\dinil{-}g(x)\in[\lambda-g(x),\alpha-g(x)]$ for all $x\in U$.
    \end{enumerate}
\end{lemma}
\begin{proof}
    We will see that $\dinil{+}g(x)\in[\lambda-g(x),\alpha-g(x)]$ for all $x\in U$ implies that $\diniu{+}g(x)\in[\lambda-g(x),\alpha-g(x)]$ for all $x\in U$; the remaining implications are similar.

    Suppose for contradiction there is some $x_0\in U$ such that $\diniu{+}g(x_0)\notin[\lambda-g(x_0),\alpha-g(x_0)]$.
    If $\diniu{+}g(x_0)<\lambda-g(x_0)$ this is immediate, so we assume $\diniu{+}g(x_0)>\alpha-g(x_0)$.
    Then there is some $\epsilon>0$ and $x_1$ such that
    \begin{equation*}
        \frac{g(x_1)-g(x_0)}{x_1-x_0}\geq \alpha-g(x_0)+\epsilon,
    \end{equation*}
    $[x_0,x_1]\subset U$, and $|g(y)-g(x_0)|<\epsilon/2$ for all $y\in[x_0,x_1]$.
    Then by \cref{c:dmvt}, there is some $y\in[x_0,x_1]$ such that
    \begin{equation*}
        \dinil{+}g(y)\geq \alpha-g(x_0)+\epsilon>\alpha-g(y)+\frac{\epsilon}{2}
    \end{equation*}
    a contradiction.
\end{proof}
\subsection{Bounding the intermediate dimensions}
In this section, we provide a general bound for the intermediate dimensions.
\begin{definition}
    Given $\theta\in[0,1]$, we say that a family of sets $\{U_i\}_{i=1}^\infty$ is a $(\delta,\theta)$-cover of $F$ if $F\subseteq\bigcup_{i=1}^\infty U_i$ and $\delta^{1/\theta}\leq\diam U_i\leq\delta$ for each $i\in\N$.
\end{definition}
We take $\delta^{1/0}=0$.
For convenience, given a cover $\mathcal{U}=\{U_i\}_{i=1}^\infty$, we define the \defn{$s$-cost} of the cover by $\mathcal{C}^s(\mathcal{U})\coloneqq\sum_{i=1}^\infty(\diam U_i)^s$.
With this terminology, we recall that for $0\leq\theta\leq 1$, the upper intermediate dimensions are given by
\begin{equation*}
    \overline{\dim}_\theta F=\inf\Bigl\{s\geq 0:\exists\delta_0>0\,\forall0<\delta\leq\delta_0\,\exists\text{$(\delta,\theta)$-cover }\mathcal{U}\text{ of }F \text{ s.t. }\mathcal{C}^s(\mathcal{U})\leq 1\Bigr\}
\end{equation*}
and the lower intermediate dimensions are given by
\begin{equation*}
    \underline{\dim}_\theta F=\inf\Bigl\{s\geq 0:\forall\delta_0>0\,\exists0<\delta\leq\delta_0\text{ and $(\delta,\theta)$-cover }\mathcal{U}\text{ of }F\text{ s.t. }\mathcal{C}^s(\mathcal{U})\leq 1\Bigr\}.
\end{equation*}
Of course, it suffices to show that $\mathcal{C}^s(\mathcal{U})\leq M$ for some constant $M$ not depending on $\delta$.

We now prove the following general bound for subsets of $\R^d$.
The proof uses a similar strategy to the proof of the bound \cite[Theorem~3.5]{ban2020} for a more general family of dimensions known as the $\Phi$-intermediate dimensions.
The idea behind the proof is as follows.
We will bound $\overline{\dim}_{\theta+\epsilon}F$ in terms of $\theta$, $\overline{\dim}_\theta F$, and the Assouad and lower dimensions of $F$.
Given an optimal cover for $[\delta^{1/\theta},\delta]$, we want to convert this to a cover for the smaller range of scales $[\delta^{\beta/(\theta+\epsilon)},\delta^\beta]\subset[\delta^{1/\theta},\delta]$.
We then use the Assouad dimension to replace the sets with large diameter with sets with smaller diameter (corresponding to the indices in $I_3$), and the lower dimension to optimally cover the sets with small diameter (corresponding to the indices in $I_1$).
The remaining elements of the cover remain essentially the same.
The parameter $\beta$ is chosen to optimize this process.

In order to obtain bounds corresponding to the lower dimension, we find it convenient to use lower dimensions of measures.
If $\mu$ is a Borel probability measure, the \defn{lower dimension} of $\mu$ is given by
\begin{align*}
    \dimL\mu = \sup\Bigl\{\lambda:&(\exists C>0)\,(\forall 0<r<R\leq\diam(\supp\mu)\text{ and }x\in \supp\mu),\\*
                                  &\frac{\mu(B(x,R))}{\mu(B(x,r))}\geq C\left(\frac{R}{r}\right)^\lambda\Bigr\}.
\end{align*}
We refer the reader to \cite[Section~4.1]{fra2021} for more details.

We recall that a measure $\mu$ is \defn{doubling} if there exists $M \geq 1$, called the \defn{doubling constant}, such that $\mu(B(x,2r)) \leq M\mu(B(x,r))$ for all $x \in\supp(\mu)$ and $r>0$.
We also recall that $\mathcal{H}(\lambda,\alpha)$ is defined in \cref{d:h-class}.
\begin{theorem}\label{p:int-dim-bound}
    Let $F\subset\R^d$ be any bounded set with $\lambda = \dimL F$, $\alpha = \dimA F$.
    Write $\overline{h}(\theta)=\overline{\dim}_\theta F$ and $\underline{h}(\theta)=\underline{\dim}_\theta F$ and let $h\in\{\underline{h},\overline{h}\}$.
    Then, $h\in\mathcal{H}(\lambda,\alpha)$.
    In particular, if $h(\theta)\in\{\lambda,\alpha\}$ for some $0<\theta\leq 1$, then $h(\theta)$ is constant on $(0,1]$.
\end{theorem}
\begin{proof}
    We prove this for $h(\theta)=\overline{\dim}_\theta F$; the case when $h(\theta)=\underline{\dim}_\theta F$ is similar.
    We also assume that $\lambda < \alpha$, or else the result is trivial.
    First let $\theta \in (0,1)$ and $\epsilon \in (0,1-\theta)$, and let $\eta,\beta$ be the unique solutions to the equations
    \begin{align*}
        \alpha - &h(\theta) - \beta (\alpha - h(\theta) - \eta) = 0 & \frac{\beta}{\theta + \epsilon}&(h(\theta) + \eta - \lambda) + \frac{\lambda - h(\theta)}{\theta} = 0.
    \end{align*}
    One can verify that $\eta$ and $\beta$ are given by
    \begin{align*}
        \eta &= \frac{( h(\theta) - \lambda) (\alpha -h(\theta))\epsilon}{(h(\theta) - \lambda)\epsilon + (\alpha - \lambda)\theta} & \beta &= \frac{(h(\theta) - \lambda)\epsilon}{(\alpha - \lambda)\theta}+1.
    \end{align*}
    Now for $s>h(\theta)$, let $s' \in (h(\theta),s)$, $\alpha' > \alpha$ and $\lambda' < \lambda$ satisfy
    \begin{align*}
        \alpha' - &s - \beta (\alpha' - s - \eta) > 0 & \frac{\beta}{\theta + \epsilon}&(s + \eta - \lambda') + \frac{\lambda' - s'}{\theta} > 0.
    \end{align*}
    For all sufficiently small $\delta \in (0,1)$ there exists a $(\delta,\theta)$-cover $\{U_i\}_{i \in I}$ of $F$ whose $s'$-cost is less than 1.
    Define
    \begin{align*}
        I_1 &= \{ i \in I : \diam U_i < \delta^{\frac{\beta}{\theta + \epsilon}} \}\\
        I_2 &= \{ i \in I :  \delta^{\frac{\beta}{\theta + \epsilon}} \leq  \diam U_i \leq \delta^\beta / 2\}\\
        I_3 &= \{ i \in I :   \diam U_i > \delta^\beta / 2\}.
    \end{align*}

    There exists $C > 0$ such that for all $0 < r \leq 2R$, any set of diameter $R$ contained in $F$ can be covered by $\lfloor C(R/r)^{\alpha'} \rfloor$ balls of diameter $r$.
    Then for $k \in I_3$, let
    \begin{equation*}
        B_{k,1}, \ldots, B_{k,\lfloor C((2\diam U_k)/\delta^{\beta})^{\alpha'} \rfloor}
    \end{equation*}
    satisfy
    \begin{equation*}
        \diam B_{k,i}=\delta^\beta\quad\text{and}\quad S_{\delta^{\beta/(\theta+\epsilon)}}(U_k)\cap F\subset\bigcup_{i=1}^{\lfloor C((2\diam U_k)/\delta^{\beta})^{\alpha'} \rfloor}B_{k,i},
    \end{equation*}
    where $S_r(U)$ denotes the $r$-neighbourhood of the set $U$.
    Let $z_1,\ldots,z_K$ be a maximal $4\delta^{\beta/(\theta + \epsilon)}$-separated subset of
    \begin{equation*}
        F \setminus \left( \bigcup_{i \in I_2 \cup I_3} S_{\delta^{\frac{\beta}{\theta + \epsilon}}}(U_i) \right).
    \end{equation*}
    Set
    \begin{align*}
        \mathcal{U}_1 &\coloneqq \{ B(z_m, 5\delta^{\frac{\beta}{\theta + \epsilon}}) : 1 \leq m \leq K \}\\
        \mathcal{U}_2 &\coloneqq \{ S_{\delta^{\frac{\beta}{\theta + \epsilon}}}(U_j) : j \in I_2 \}\\
        \mathcal{U}_3 &\coloneqq \bigcup_{k\in I_3}\bigl\{B_{k,\ell}:\ell=1,\ldots,\lfloor C((2\diam U_k)/\delta^{\beta})^{\alpha'} \rfloor\bigr\}.
    \end{align*}
    Then for sufficiently small $\delta$,
    \begin{equation}\label{e:general-cover}
        \mathcal{U}\coloneqq\mathcal{U}_1\cup\mathcal{U}_2\cup\mathcal{U}_3
    \end{equation}
    is a $(\delta^\beta , \theta + \epsilon)$-cover of $F$.

    We now bound the $(s+\eta)$-cost of $\mathcal{U}$ independently of $\delta$.
    First consider $\mathcal{U}_1$.
    By \cite[Theorem~2]{bg2000}, there exists a doubling Borel probability measure $\mu$ with $\supp\mu=\overline{F}$ and $\dimL \mu  \in (\lambda',\lambda]$.
    Let $M$ be a doubling constant for $\mu$.
    In particular, there is $c > 0$ such that if $0 < r < R \leq \diam F$ and $x \in F$ then
    \begin{equation*}
        \frac{\mu(B(x,R))}{\mu(B(x,r))} \geq c \left(\frac{R}{r}\right)^{\lambda'}.
    \end{equation*}
    For $m \in \{1,\ldots,K\}$ let
    \begin{equation*}
        J_m\coloneqq \{i\in I_1:U_i \cap B(z_m,\delta^{\beta/(\theta + \epsilon)}) \neq \varnothing\}.
    \end{equation*}
    If $i \in J_m$, fixing $x_{i,m} \in U_i \cap B(z_m,\delta^{\beta/(\theta + \epsilon)})$,
    \begin{align*}
        \mu(U_i) &\leq \mu(B(x_{i,m},2\diam(U_i))) \leq c^{-1} \mu(B(x_{i,m},2\delta^{\frac{\beta}{\theta + \epsilon}})) \left(\frac{ \delta^{\frac{\beta}{\theta + \epsilon}} }{\diam U_i}\right)^{-\lambda'} \\
        &\leq c^{-1} \mu(B(z_m,4\delta^{\frac{\beta}{\theta + \epsilon}})) \left(\frac{ \delta^{\frac{\beta}{\theta + \epsilon}} }{\diam U_i}\right)^{-\lambda'} \leq c^{-1} M^2\mu(B(z_m,\delta^{\frac{\beta}{\theta + \epsilon}})) \left(\frac{ \delta^{\frac{\beta}{\theta + \epsilon}} }{\diam U_i}\right)^{-\lambda'}.
    \end{align*}
    Then
    \begin{equation*}
        \mu(B(z_m,\delta^{\frac{\beta}{\theta + \epsilon}}))  \leq \sum_{i\in J_m} \mu(U_i) \leq c^{-1} M^2 \mu(B(z_m,\delta^{\frac{\beta}{\theta + \epsilon}})) \delta^{\frac{-\lambda' \beta}{\theta + \epsilon}} \sum_{i\in J_m} (\diam U_i)^{\lambda'}.
    \end{equation*}
    Note that $\mu(B(z_m,\delta^{\frac{\beta}{\theta + \epsilon}})) > 0$ since $\supp\mu=\overline{F}$.
    Moreover, if $i \in I_1$, then $U_i$ intersects at most one of the balls of radius $\delta^{\beta/(\theta + \epsilon)}$, so for sufficiently small $\delta$,
    \begin{align*}
        \mathcal{C}^{s+\eta}(\mathcal{U}_1)&=\sum_{i \in I_1} \diam(B(z_i, 5\delta^{\frac{\beta}{\theta + \epsilon}}))^{s + \eta}\\
                                           &= \sum_{m = 1}^K (10\delta^{\frac{\beta}{\theta + \epsilon}})^{s + \eta}\leq 10^{s + \eta} c^{-1} M^2 \delta^{\frac{\beta}{\theta + \epsilon}(s + \eta - \lambda')} \sum_{i \in I} (\diam U_i)^{\lambda'} \\
                                                                                   &\leq 10^{s + \eta} c^{-1} M^2 \delta^{\frac{\beta}{\theta + \epsilon}(s + \eta - \lambda')} \delta^{\frac{\lambda' - s'}{\theta}} \sum_{i \in I} (\diam U_i)^{s'} \leq 10^{s + \eta} c^{-1} M^2.
    \end{align*}

    Next, consider $\mathcal{U}_2$:
    \begin{equation*}
        \mathcal{C}^{s+\eta}(\mathcal{U}_2)=\sum_{j \in I_2} \diam(S_{\delta^{\frac{\beta}{\theta + \epsilon}}}(U_j))^{s+\eta} \leq  \sum_{j \in I_2} (3 \diam U_i)^{s + \eta} \leq 3^{s + \eta}.
    \end{equation*}

    Finally, consider $\mathcal{U}_3$.
    Since $\diam U_k\leq\delta$,
    \begin{align*}
        \mathcal{C}^{s+\eta}(\mathcal{U}_3)&=\sum_{k \in I_3} \sum_{\ell=1}^{\lfloor C((2\diam U_k)/\delta^{\beta})^{\alpha'} \rfloor} (\diam B_{k,\ell})^{s+\eta}\leq \sum_{k \in I_3} 2^{\alpha'} C (\diam U_k)^{\alpha'} \delta^{-\beta \alpha'} \delta^{\beta (s+\eta)}\\
                                           &\leq 2^{\alpha'} C \sum_{k \in I_3} (\diam U_k)^s \delta^{\alpha' - s - \beta \alpha' + \beta(s+\eta)}\leq 2^{\alpha'} C \sum_{k \in I} (\diam U_k)^s \leq 2^{\alpha'} C.
    \end{align*}
    Thus $\mathcal{C}^{s+\eta}(\mathcal{U})\leq 10^{s + \eta} c^{-1} M^2 + 3^{s + \eta} + 2^{\alpha'} C$ which does not depend on $\delta$.

    Since $s > h$ was arbitrary, we have shown that $h(\theta + \epsilon) \leq h(\theta) + \eta$.
    Dividing through by $\epsilon$ gives
    \begin{equation}\label{e:h-epsilon-bound}
        \frac{h(\theta+\epsilon)-h(\theta)}{\epsilon}\leq \frac{ ( h(\theta) - \lambda) (\alpha -h(\theta))}{(h(\theta) - \lambda)\epsilon + (\alpha - \lambda)\theta} .
    \end{equation}
    Passing to the limit, we verify \cref{e:h-bound-general}.
    That $h$ is non-decreasing follows immediately from the definition.

    Moreover, \cref{e:h-epsilon-bound} implies that $h$ is continuous on $(0,1)$.
    To see that $h(\theta)$ is continuous at 1, with minor modifications to the above proof we can take $\epsilon=1-\theta$.

    The particular cases $h(\theta)=\lambda$ or $h(\theta)=\alpha$ for some $\theta\in(0,1]$ follow directly from \cref{e:h-epsilon-bound}.
\end{proof}
\begin{remark}
    We also observe that the bound \cref{e:h-epsilon-bound} for any $\epsilon\in[0,1-\theta]$ can be obtained directly by solving the differential equation corresponding to \cref{e:h-bound-general} and applying \cref{l:c1-bound} (using continuity of $h(\theta)$ on $(0,1]$).
\end{remark}
As an application, we can get a general lower bound for $h(\theta)$ in terms of the lower, box, and Assouad dimensions of the set.
Recall that $\overline{\dim}_1 F=\dimuB F$ and $\underline{\dim}_1 F=\dimlB F$.
\begin{corollary}\label{c:lower-bound-general}
    Let $\dimL F=\lambda$ and $\dimA F=\alpha$.
    For $h(\theta)=\overline{\dim}_\theta F$ or $h(\theta)=\underline{\dim}_\theta F$,
    \begin{equation*}
        h(\theta)\geq\frac{\alpha \theta(h(1)-\lambda)+\lambda(\alpha-h(1))}{\theta(h(1)-\lambda)+(\alpha-h(1))}
    \end{equation*}
    for all $\theta\in(0,1]$.
\end{corollary}
\begin{proof}
    Substitute $\epsilon=1-\theta$ into \cref{e:h-epsilon-bound}.
    Rearranging, we obtain the desired result.
\end{proof}
\begin{remark}
    Differentiating twice, we see that the general lower bound in \cref{c:lower-bound-general} is a concave function of $\theta$.
    As $\dimB F$ approaches $\dimA F$ (resp.~$\dimL F$), the lower bound converges pointwise to the Assouad (resp.~lower) dimension for all $\theta>0$.
    The lower bound is always equal to $\dimL F$ for $\theta=0$.
    Plots of the lower bound for particular parameters are given in \cref{f:lower-bound-curves}.
    The bounds \cref{e:h-epsilon-bound} and \cref{c:lower-bound-general} have subsequently been proved in a class of metric spaces which are more general than $\R^d$, see \cite[Theorem~3.12 and Corollary~3.14]{ban2020}.
\end{remark}
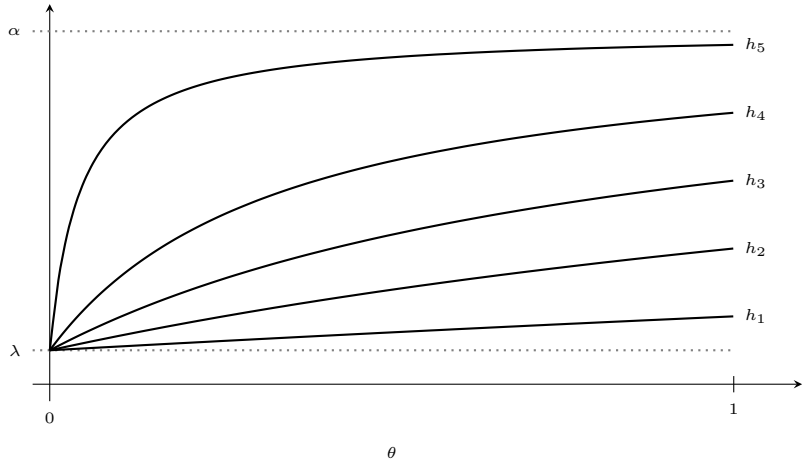
\begin{figure}[t]
    \centering
    \def\lowerdim{0.05}
\def\assouad{0.52}
\begin{tikzpicture}[>=stealth,scale=9,font=\tiny]
    \draw[->] (-0.025,0) -- (1.1,0);
    \draw[->] (0,-0.025) node[below]{$0$} -- (0,.56);
    \draw[] (1,-0.012) node[below]{$1$} -- (1,0.012);
    
    \node (xlabel) at (0.5,-0.1) {$\theta$};

    \draw[thick,dotted,gray] (-0.025,\lowerdim) node[left,black]{$\lambda$} -- (1.0,\lowerdim);
    \draw[thick,dotted,gray] (-0.025,\assouad) node[left,black]{$\alpha$} -- (1.0,\assouad);

    \foreach \idx/\boxdim in {1/0.1,2/0.2,3/0.3,4/0.4,5/0.5} {
        \draw[domain=0:1, smooth, variable=\t, samples=80,thick] plot ({\t}, {(\assouad*\t*(\boxdim-\lowerdim)+\lowerdim*(\assouad-\boxdim))/(\t*(\boxdim-\lowerdim)+(\assouad-\boxdim))}) node[black,right]{$h_{\idx}$};
    }
\end{tikzpicture}
    \caption{Plots of the bound in \cref{c:lower-bound-general} for $\lambda=0.05$, $\alpha=0.52$, and box dimensions $h_i(1)=i/10$ for $i=1,2,\ldots,5$.}
    \label{f:lower-bound-curves}
\end{figure}
\subsection{Intermediate dimensions of homogeneous Moran sets}
We first define the main object used in our construction, which we will call homogeneous Moran sets.
The construction is analogous to the usual $2^d$-corner Cantor set, except that the subdivision ratios need not be the same at each level.

Fix $\mathcal{I}=\{0,1\}^d$.
We write $\mathcal{I}^*=\bigcup_{n=0}^\infty\mathcal{I}^n$, and denote the word of length $0$ by $\varnothing$.
Suppose we are given a sequence $\bm{r}=(r_n)_{n=1}^\infty$ with $0<r_n\leq 1/2$ for each $n\in\N$.
Then for each $n$ and $\bm{i}\in\mathcal{I}$, we define $S^n_{\bm{i}}\colon\R^d\to\R^d$ by
\begin{equation*}
    S^n_{\bm{i}}(x)\coloneqq r_n x+b^n_{\bm{i}}
\end{equation*}
where $b^n_{\bm{i}}\in\R^d$ has
\begin{equation*}
    (b^n_{\bm{i}})^{(j)} =
    \begin{cases}
        0 &: \bm{i}^{(j)}=0\\*
        1-r_n &: \bm{i}^{(j)}=1
    \end{cases}.
\end{equation*}
Given $\sigma=(\bm{i}_1,\ldots,\bm{i}_n)\in\mathcal{I}^n$, we write $S_\sigma=S^1_{\bm{i}_1}\circ\cdots\circ S^n_{\bm{i}_n}$.
Then set
\begin{equation*}
    C_n=\bigcup_{\sigma\in\mathcal{I}^n}S_\sigma([0,1]^d)\qquad\text{and}\qquad C=C(\bm{r})\coloneqq\bigcap_{n=1}^\infty C_n.
\end{equation*}
We refer to the set $C$ as a \defn{homogeneous Moran set}.
Note that $C_n$ consists of $2^{dn}$ hypercubes each with diameter $\rho_n\coloneqq r_1\cdots r_n$ (with respect to the max norm).

Given $\delta>0$, let $k=k(\delta)$ be such that $\rho_k\leq \delta < \rho_{k-1}$.
We then define
\begin{equation*}
    s(\delta)=s_{\bm{r}}(\delta)\coloneqq\frac{k(\delta)\cdot d\log 2}{-\log\delta}.
\end{equation*}
One can interpret $s(\delta)$ as the best candidate for the ``box dimension at scale $\delta$''.

We now prove the following key covering lemma for intermediate dimensions.
This result essentially shows that the optimal covers for a homogeneous Moran set can be taken to consist of balls all of the same diameter.
\begin{lemma}\label{l:flat-covers}
    Let $\theta\in(0,1]$ be arbitrary.
    Then for all $\delta>0$ sufficiently small, with $t=\inf_{\phi\in[\delta^{1/\theta},\delta]}s(\phi)$,
    \begin{equation*}
        4^{-d}\leq \inf\Bigl\{\mathcal{C}^t(\mathcal{U}):\mathcal{U}\text{ is a $(\delta,\theta)$-cover of $C$}\Bigr\}\leq 1.
    \end{equation*}
\end{lemma}
\begin{proof}
    Let $\mu$ denote the uniform Bernoulli measure on $C$.
    Let $U$ be a set with $\delta^{1/\theta}\leq\diam U\leq\delta$, and let $k$ be such that $\rho_{k}\leq\diam U<\rho_{k-1}$.
    Note that $(\diam U)^{s(\diam U)}=2^{-kd}$.
    Then since $U$ intersects at most $4^d$ hypercubes in $C_k$,
    \begin{equation*}
        \mu(U)\leq 4^d\cdot 2^{-kd}= 4^d\cdot (\diam U)^{s(\diam U)}\leq 4^d(\diam U)^t.
    \end{equation*}
    In particular, if $\mathcal{U}$ is an arbitrary $(\delta,\theta)$-cover of $C$,
    \begin{equation*}
        1=\mu(C)\leq\sum_{U\in\mathcal{U}}\mu(U)\leq 4^d\sum_{U\in\mathcal{U}}(\diam U)^t
    \end{equation*}
    so that $\mathcal{C}^t(\mathcal{U})\geq 4^{-d}$.

    Conversely, since $s(\delta)$ is continuous and increasing on each interval $[\rho_k,\rho_{k-1})$, there is $\phi\in[\delta^{1/\theta},\delta]$ such that $t=s(\phi)=\frac{k(\phi)\cdot d\log 2}{-\log\phi}$.
    For each $y=(j_1,\ldots,j_d)\in\{0,1\}^d$ and $\sigma\in\mathcal{I}^*$, let $E_{\sigma,\phi}(y)$ denote the hypercube with side length $\phi$ contained in $S_\sigma([0,1]^d)$, with edges aligned with the coordinate axes, and containing the point $S_\sigma(y)$.
    Since $\phi\geq\rho_{k(\phi)}$,
    \begin{equation*}
        \mathcal{V}\coloneqq\bigcup_{\sigma\in\mathcal{I}^{k(\phi)-1}}\{E_{\sigma,\phi}(y):y\in\{0,1\}^d\}
    \end{equation*}
    is a cover for $C_{k(\phi)}$, and therefore $C$, consisting of $2^{k(\phi)\cdot d}$ hypercubes each with diameter $\phi$.
    Thus $\mathcal{C}^t(\mathcal{V})=1$.
\end{proof}
As a direct application, we have the following formula for the intermediate dimensions of $C$.
\begin{proposition}\label{p:int-dim}
    For any $\theta\in(0,1]$,
    \begin{equation*}
        \overline{\dim}_\theta C = \limsup_{\delta\to 0}\bigl(\inf_{\phi\in[\delta^{1/\theta},\delta]}s(\phi)\bigr)
    \end{equation*}
    and
    \begin{equation*}
        \dimH C = \underline{\dim}_\theta C = \dimlB C=\liminf_{\delta\to 0}s(\delta).
    \end{equation*}
\end{proposition}
\section{Constructions with Moran sets}
In this section, we prove the converse direction to \cref{it:general-form}.

In \cref{ss:Moran-cnst}, we will establish a general strategy for constructing homogeneous Moran sets.
We first introduce the following definition, which is in some sense analogous to the definition of $\mathcal{H}(\lambda,\alpha)$.
\begin{definition}\label{d:g-class}
    Given $0\leq\lambda\leq\alpha\leq d$, we write $\mathcal{G}(\lambda,\alpha)$ to denote the functions $g\colon\R\to[\lambda,\alpha]$ which are continuous and satisfy
    \begin{equation}
        \diniu{+}g(x)\in[\lambda-g(x),\alpha-g(x)]
    \end{equation}
    for all $x\in\R$.
\end{definition}
We will essentially show that for any function $g\in\mathcal{G}(0,d)$, there exists a homogeneous Moran set such that $s(\delta)\approx g(\log\log(1/\delta))$.
The transformation $\delta\mapsto\log\log(1/\delta)$ is useful since it converts the exponentiation map $\delta\mapsto\delta^{1/\theta}$ into addition $x\mapsto x+\log(1/\theta)$.

In order to construct such a set, it suffices to define the corresponding contraction ratios by ``discretizing'' the function $g$.
In particular, in \cref{l:exact-construction}, we show that there exists a sequence of contractions $\bm{r}$ such that the corresponding covering numbers $s_{\bm{r}}(\delta)$ are close to $g(\log\log(1/\delta))$ in the precise sense given in \cref{e:g-disc}.
Of course, depending on the choice of the function $g$, this bound may be impossible to attain for small $x$.
Thus we begin with a function $\tilde{g}$ and then translate it by some constant amount.
The contraction ratios are then used to define a corresponding Moran set $C$, and \cref{e:g-disc} is useful to prove dimension results for the Moran set $C$.

Then, in \cref{ss:pres-upper} and \cref{ss:pres}, we use this technique to construct Moran sets with the desired properties.
In \cref{t:upper-h-form}, we construct the function $g$ depending on some $h\in\mathcal{H}(\lambda,\alpha)$ such that the corresponding Moran set has the desired dimension formulas.
This construction is also used in \cref{t:gen-h-form}, where we use the sequence of contraction ratios provided by \cref{l:exact-construction} directly.
Here, translations of the function $g$ are used to define an inhomogeneous Moran set which ``locally'' looks like the Moran set $C$, but with a much greater amount of uniformity between scales (so that the intermediate dimensions exist).
Finally, these results are combined in \cref{c:upper-lower-match} to obtain a proof of \cref{it:general-form}.

\subsection{Constructing homogeneous Moran sets}\label{ss:Moran-cnst}
We first describe a general strategy to construct homogeneous Moran sets.
\begin{lemma}\label{l:scale-bound}
    Let $0\leq\lambda\leq\alpha\leq d$ and let $g \colon \R \to [\lambda,\alpha]$.
    Then $g\in\mathcal{G}(\lambda,\alpha)$ if and only if for all $x_0\in\R$ and $x>0$,
    \begin{equation*}
        \lambda-(\lambda-g(x_0))\exp(-x)\leq g(x_0+x)\leq \alpha-(\alpha-g(x_0))\exp(-x).
    \end{equation*}
\end{lemma}
\begin{proof}
    This is a direct application of \cref{l:c1-bound}.
\end{proof}
\begin{definition}
    Given a sequence of functions $(f_k)_{k=1}^\infty$ each defined on some interval $[0,a_k]$, the \defn{concatenation} of $(f_k)_{k=1}^\infty$ is the function $f\colon(-\infty,\sum_{k=1}^\infty a_k)\to\R$ given as follows: for each $x>0$ with $\sum_{j=0}^{k-1} a_j<x\leq\sum_{j=0}^{k}a_j$ where $a_0=0$, we define
    \begin{equation*}
        f(x)=f_k\left(x-\sum_{j=0}^{k-1} a_j\right)
    \end{equation*}
    and for $x\leq 0$ we define $f(x)=f_1(0)$.
\end{definition}

Given a function $g\in\mathcal{G}(\lambda,\alpha)$ and $w\in\R$, we define the \emph{offset} $\kappa_w(g)\in\mathcal{G}(\lambda,\alpha)$ by
\begin{equation*}
    \kappa_w(g)(x)=\begin{cases}
        g(x-w) &: x\geq w,\\
        g(0) &: x\leq w.
    \end{cases}
\end{equation*}
We also say that a function $g\in\mathcal{G}(\lambda,\alpha)$ is \emph{rapidly decreasing} if there is a $y\in\R$ and a constant $C>0$ so that for all $x\geq y$,
\begin{equation}\label{e:rapidly-decreasing}
    g(x)\leq g(y)\exp(y-x)+C\exp(-x).
\end{equation}
Note that if $g$ is rapidly decreasing, then $\lim_{x\to\infty}g(x)=0$.
Moreover, for all $w\in\R$, $g$ is not rapidly decreasing if and only if $\kappa_w(g)$ is not rapidly decreasing.

The following lemma is stated to be useful in the proof of \cref{t:gen-h-form}, where many offsets of the same function will be required.
\begin{lemma}\label{l:exact-construction}
    Let $0\leq\lambda<\alpha\leq d$ and let $\tilde{g}\in\mathcal{G}(\lambda,\alpha)$.
    Suppose $\tilde{g}$ is not rapidly decreasing.
    Then there is a constant $w_0\in\R$ depending only on $\tilde{g}(0)$ and $d$ such that for all $w\geq w_0$, there exists a sequence $\bm{r}\coloneqq (r_j)_{j=1}^\infty\subset(0,1/2]$ so that $g\coloneqq\kappa_{w}(\tilde{g})$ satisfies
    \begin{equation}\label{e:g-disc}
        |s_{\bm{r}}(\exp(-\exp(x)))-g(x)|\leq d\log(2)\cdot\exp(-x)
    \end{equation}
    for all $x\geq w_0$.
\end{lemma}
\begin{proof}
    Noting that $\tilde{g}(0)\in(0,d)$, choose $r_1$ such that $\frac{2d\log(2)}{\log(1/r_1)}=\tilde{g}(0)$.
    Then let $w_0=\log\log(1/r_1)$, let $w\geq w_0$ be arbitrary, and let $g=\kappa_w(\tilde{g})$.
    Since $\tilde{g}$ is not rapidly decreasing, $g$ is also not rapidly decreasing so by \cref{e:rapidly-decreasing} for every $y\in\R$ there is a minimal $\psi(y)>y$ so that
    \begin{equation*}
        g(y)\exp(y-\psi(y))=g(\psi(y))-d\log(2)\cdot\exp(-\psi(y)).
    \end{equation*}

    Now set $x_1=w_0$ and, inductively, set $x_{k+1}=\psi(x_k)$ for each $k\in\N$.
    Let $\rho_k=\exp(-\exp(x_k))$ denote the corresponding scales (note that $\rho_1 = r_1$), and set $r_k\coloneqq\rho_k/\rho_{k-1}$ for $k\geq 2$.
    Observe that $r_k\in(0,1)$ for all $k$.
    Thus for $0<\delta\leq r_1$, if $k$ is such that $\rho_k<\delta\leq\rho_{k-1}$, we set
    \begin{equation*}
        \overline{s}(\delta)=\frac{kd\log 2}{\log(1/\delta)}.
    \end{equation*}

    We will prove by induction that for each $k\in\N$ we have $r_k\in(0,1/2]$, $\overline{s}(\rho_k)=g(x_k)$, and
    \begin{equation}\label{e:bound}
        g(x)-d\log(2)\exp(-x)\leq\overline{s}(\exp(-\exp(x)))\leq g(x)
    \end{equation}
    holds for all $x\in[x_1,x_k]$.
    From this, the result follows.

    We first note that, by construction, $r_1\in(0,1/2]$ and $\overline{s}(\rho_1)=g(x_1)=\tilde{g}(0)$.
    In general, suppose the hypothesis holds for $k\in\N$.
    By definition of $\psi$ and the fact that $g(x_k)=\overline{s}(\rho_k)$,
    \begin{align*}
        g(x_{k+1}) &=\overline{s}(\rho_k)\exp(-x_{k+1}+x_k)+d\log(2)\exp(-x_{k+1})\\
                   &=\frac{d(k+1)\log 2}{\exp(x_k)}\cdot\exp(-x_{k+1})\exp(x_k)+d\log(2)\exp(-x_{k+1})\\
                   &= \frac{d(k+2)\log 2}{\exp(x_{k+1})}= \overline{s}(\rho_{k+1}).
    \end{align*}
    Moreover, by \cref{l:scale-bound}, $g(x)\geq g(x_k)\exp(-x+x_k)$ for all $x\geq x_k$ so that \cref{e:bound} follows for $x\in[x_k,x_{k+1}]$ by the minimality of $x_{k+1}$ in the definition of $\psi$.
    Finally, again by \cref{l:scale-bound}, $g(x_{k+1})\leq d-(d-g(x_k))\exp(-{x_{k+1}}+x_k)$.
    Substituting, this implies that
    \begin{equation*}
        \frac{d(k+2)\log 2}{\log(1/\rho_{k+1})}\leq d-\left(d-\frac{d(k+1)\log 2}{\log(1/\rho_k)}\right)\cdot\frac{\log(1/\rho_k)}{\log(1/\rho_{k+1})}
    \end{equation*}
    which after simplification gives that $\rho_{k+1}\leq\rho_k/2$, i.e. $r_{k+1}\leq 1/2$.
\end{proof}
\begin{remark}
    If instead $g$ is rapidly decreasing, then the function $g$ is decays faster than any function $s_{\bm{r}}$ for a sequence $\bm{r}\subset(0,1/2]$.
\end{remark}
\begin{remark}
    The bound \cref{e:bound} is optimal since $s(\delta)$ has discontinuities of size $\frac{d\log 2}{\log(1/\delta)}$.
\end{remark}
We now use the sequence $\bm{r}$ constructed in the previous lemma to define a homogeneous Moran set $C$, and prove that it satisfies the correct properties.
Recall that $\mathcal{G}$ is defined in \cref{d:g-class}.
\begin{lemma}\label{l:Moran-formula}
    Let $g\in\mathcal{G}(0,d)$ and suppose $\bm{r}=(r_j)_{j=1}^\infty\subset(0,1/2]$ satisfies
    \begin{equation}\label{e:tight-bound}
        |s_{\bm{r}}(\exp(-\exp(x)))-g(x)|\leq d\log(2)\cdot\exp(-x)
    \end{equation}
    for all $x$ sufficiently large.
    Then the corresponding homogeneous Moran set $C=C(\bm{r})\subset\R^d$ satisfies:
    \begin{enumerate}[nl,r]
        \item\label{im:int-dim} $\displaystyle\overline{\dim}_\theta C =\limsup_{x\to\infty}\left(\inf_{y\in[x,x+\log(1/\theta)]}g(y)\right)$ for $\theta\in(0,1]$,
        \item\label{im:h-dim} $\displaystyle\underline{\dim}_\theta C=\dimH C=\liminf_{x\to\infty}g(x)$ for $\theta\in(0,1]$,
        \item\label{im:A-dim-gen} $\dimA C\leq\limsup_{x\to\infty}(\diniu{+}g(x)+g(x))$, and
        \item\label{im:L-dim-gen} $\dimL C\geq\liminf_{x\to\infty}(\diniu{+}g(x)+g(x))$.
    \end{enumerate}
    Moreover, suppose $\psi \colon \R \to \R^+$ is any function such that $\lim_{x\to\infty}(\exp(\psi(x))-\exp(x))=\infty$.
    Then
    \begin{enumerate}[nl,r,resume]
        \item\label{im:A-dim-block} $\displaystyle\dimA C\geq\limsup_{x\to\infty}\left(\inf_{y\in[x,\psi(x)]}(\diniu{+}g(y)+g(y))\right)$, and
        \item\label{im:L-dim-block} $\displaystyle\dimL C\leq\liminf_{x\to\infty}\left(\sup_{y\in[x,\psi(x)]}(\diniu{+}g(y)+g(y))\right)$.
    \end{enumerate}
\end{lemma}
\begin{proof}
    We first observe that \cref{im:int-dim} and \cref{im:h-dim} follow immediately from \cref{p:int-dim}.
    We will verify \cref{im:A-dim-gen} and \cref{im:A-dim-block}; \cref{im:L-dim-gen} and \cref{im:L-dim-block} are given by an analogous argument.

    We first establish a general formula for the Assouad dimension of $C$ in terms of the numbers $s_{\bm{r}}(\delta)$.
    Suppose $0<\delta_1\leq\delta_2$ are arbitrary.
    Then the number of subdivision steps between scales $\delta_1$ and $\delta_2$, up to an error of size $2$, is
    \begin{equation*}
        \frac{s_{\bm{r}}(\delta_1)\log(1/\delta_1)-s_{\bm{r}}(\delta_2)\log(1/\delta_2)}{d\log 2}.
    \end{equation*}
    Thus there is a bounded function $h(x,\delta_1,\delta_2)$ such that
    \begin{equation*}
        \frac{\log N_{\delta_1}(B(x,\delta_2)\cap C)}{\log(\delta_2/\delta_1)}=\frac{s_{\bm{r}}(\delta_1)\log(1/\delta_1)-s_{\bm{r}}(\delta_2)\log(1/\delta_2)+h(x,\delta_1,\delta_2)}{\log(1/\delta_1)-\log(1/\delta_2)}.
    \end{equation*}
    Therefore by the definition of the Assouad dimension, for all $\delta_0\in(0,1)$,
    \begin{equation}\label{e:assouadlim}
        \dimA C=\lim_{\epsilon\to 0}\sup_{\substack{0<\delta_1<\delta_2<\delta_0\\\delta_1\leq\epsilon\delta_2}}\frac{s_{\bm{r}}(\delta_1)\log(1/\delta_1)-s_{\bm{r}}(\delta_2)\log(1/\delta_2)}{\log(1/\delta_1)-\log(1/\delta_2)}.
    \end{equation}
    Now we may inductively choose sequences of positive numbers $(\delta_{1,n})_{n=1}^\infty$ and $(\delta_{2,n})_{n=1}^\infty$ such that $\delta_{1,n}/\delta_{2,n}$ and $\delta_{2,n}$ converge to $0$, and
    \begin{equation*}
        \frac{s_{\bm{r}}(\delta_{1,n})\log(1/\delta_{1,n})-s_{\bm{r}}(\delta_{2,n})\log(1/\delta_{2,n})}{\log(1/\delta_{1,n})-\log(1/\delta_{2,n})} \in \Big(\dimA C - \frac{1}{n},\dimA C + \frac{1}{n}\Big)
    \end{equation*}
    for all $n \in \N$.

    For $0<x<y$, let
    \begin{equation*}
        \Phi(x,y)\coloneqq \frac{s_{\bm{r}}(\exp(-\exp(y)))-s_{\bm{r}}(\exp(-\exp(x)))}{1-\exp(x-y)}+s_{\bm{r}}(\exp(-\exp(x))).
    \end{equation*}
    Next, write $x_n=\log\log(1/\delta_{2,n})$ and $y_n=\log\log(1/\delta_{1,n})$.
    Let $\mathcal{W}$ denote the family of functions $\psi\colon\R\to\R^+$ such that $\lim_{x\to\infty}(\exp(\psi(x))-\exp(x))=\infty$.
    The condition that $\delta_{1,n}/\delta_{2,n}$ converges to $0$ is equivalent to $\exp(y_n)-\exp(x_n)$ diverging to infinity.
    Thus we may choose a function $\psi_0 \in\mathcal{W}$ so that $\psi_0(x_n)=y_n$ for infinitely many $n$. Then with some rearrangement using \cref{e:assouadlim} and the definition of $\Phi$,
    \begin{equation*}
        \limsup_{x\to\infty}\Phi(x,\psi_0(x)) \geq \dimA C.
    \end{equation*}
    Conversely, if $\psi \in\mathcal{W}$ is arbitrary, applying the substitutions $\delta_2=\exp(-\exp(x))$ and $\delta_1=\exp(-\exp(\psi(x)))$ and using the fact that $\delta_1/\delta_2$ converges to $0$ as $x \to \infty$ gives
    \begin{equation*}
        \limsup_{x\to\infty}\Phi(x,\psi(x))\leq\dimA C.
    \end{equation*}
    Therefore
    \begin{equation}\label{e:assouad-formula}
        \dimA C=\sup_{\psi\in\mathcal{W}}\limsup_{x\to\infty}\Phi(x,\psi(x)),
    \end{equation}
    and moreover the supremum is attained.

    To conclude the preliminaries, we also note, for $0<x<y$,
    \begin{equation}\label{e:error-bound}
        \frac{\exp(-y)+\exp(-x)}{1-\exp(x-y)}+\exp(-x)=\frac{2}{\exp(y)-\exp(x)}+2\exp(-x).
    \end{equation}
    This bound will be used to control the error resulting from \cref{e:tight-bound}.

    We now prove \cref{im:A-dim-gen}.
    Write
    \begin{equation*}
        \overline{\alpha} \coloneqq \limsup_{x\to\infty}(\diniu{+}g(x)+g(x)),
    \end{equation*}
    and let $\epsilon > 0$.
    Then there exists $M_\epsilon>0$ such that for all $x \geq M_\epsilon$ we have $\diniu{+}g(x)+g(x) \leq\overline{\alpha} + \epsilon$.
    For $x\geq M_\epsilon$, define $\overline{g}_x\colon [x,\infty) \to \R$ by
    \begin{equation*}
        \overline{g}_x(y) \coloneqq \overline{\alpha} + \epsilon - (\overline{\alpha} + \epsilon - g(x))\exp(x-y).
    \end{equation*}
    Then $g(x) = \overline{g}_x(x)$, and
    \begin{equation*}
        \overline{g}_x'(y) + \overline{g}_x(y) = \overline{\alpha} + \epsilon \geq \diniu{+}g(y)+g(y)
    \end{equation*}
    for all $y>x$.
    It follows from \cref{l:c1-bound} that $g(y) \leq \overline{g}_x(y)$ for all $y \geq x$.
    Now taking a function $\psi_0$ which attains the supremum in \cref{e:assouad-formula}, for all $x\geq M_\epsilon$, using \cref{e:tight-bound} and \cref{e:error-bound} combined with the condition on $\psi_0$,
    \begin{align*}
        \Phi(x,\psi_0(x)) \leq{}& \frac{g(\psi_0(x)) - g(x)}{1-\exp(x-\psi_0(x))} + g(x)\\
              &+2d\log(2)\left(\frac{1}{\exp(\psi_0(x))-\exp(x)}+\exp(-x)\right).
    \end{align*}
    Moreover, since $\psi_0\in\mathcal{W}$,
    \begin{equation*}
        \limsup_{x\to\infty}2d\log(2)\left(\frac{1}{\exp(\psi_0(x))-\exp(x)}+\exp(-x)\right)=0.
    \end{equation*}
    Thus
    \begin{align*}
        \limsup_{x\to\infty}\Phi(x,\psi_0(x)) &\leq\limsup_{x\to\infty}\left(\frac{g(\psi_0(x)) - g(x)}{1-\exp(x-\psi_0(x))} + g(x)\right)\\
                                           &\leq\limsup_{x\to\infty}\left(\frac{\overline{g}_x(\psi_0(x)) - g(x)}{1-\exp(x-\psi_0(x))} + g(x)\right)
                                           =\overline{\alpha}+\epsilon.
    \end{align*}
    But $\epsilon>0$ was arbitrary, giving the claim.

    Finally, we prove \cref{im:A-dim-block}.
    Fix any $\psi \in \mathcal{W}$ and $\epsilon>0$, and write
    \begin{equation*}
        \underline{\alpha}=\limsup_{x\to\infty}\left(\inf_{y\in[x,\psi(x)]}(\diniu{+}g(y)+g(y))\right).
    \end{equation*}
    Get a sequence $(x_k)_{k=1}^\infty$ diverging to infinity such that for all $k\in\N$,
    \begin{equation*}
        \inf_{y\in[x_k,\psi(x_k)]}(\diniu+ g(y)+g(y))\geq \underline{\alpha}-\epsilon.
    \end{equation*}
    Exactly as above, define $\underline{g}_k\colon[x_k,\infty)\to\R$ by
    \begin{equation*}
        \underline{g}_k(x) \coloneqq \underline{\alpha} + \epsilon - (\underline{\alpha} + \epsilon - g(x_k))\exp(x_k-x).
    \end{equation*}
    Then $g(x_k)=\underline{g}_k(x_k)$ and $g(x)\geq\underline{g}_k(x)$ for all $x\in[x_k,\psi(x_k)]$.
    Thus the same computations as before yield that
    \begin{align*}
        \limsup_{x\to\infty}\Phi(x,\psi_0(x))&\geq\limsup_{k\to\infty}\Phi(x_k,\psi(x_k))\\
                                          &\geq\underline{\alpha}-\epsilon-\frac{2d\log 2}{\exp(\psi(x_k))-\exp(x_k)}-2d\log(2)\exp(-x_k).
    \end{align*}
    Since $\exp(\psi(x_k))-\exp(x_k)$ diverges to infinity and $\epsilon>0$ was arbitrary, the claimed inequality follows.
\end{proof}
\begin{remark}
    In general, \cref{im:A-dim-gen} and \cref{im:L-dim-gen} will not be equalities since one would require more robust regularity assumptions about the function $g$.
\end{remark}
\subsection{Prescribing the upper intermediate dimensions}\label{ss:pres-upper}
Now, using the general construction in the previous section, we show how to construct homogeneous Moran sets with upper intermediate dimensions given by a function $h\colon[0,1]\to(0,d)$.
The main idea is to construct certain functions which have the property that for any $h(\theta)$ there are exactly two points $\{x,x+\log(1/\theta)\}$ which have value $h(\theta)$.
This ensures that the limit supremum of infima over windows $[x,x+\log(1/\theta)]$ is exactly $h(\theta)$.
See \cref{f:mountain-ctr} for a depiction of this construction.
\begin{theorem}\label{t:upper-h-form}
    Let $0\leq\lambda\leq\alpha\leq d$ and let $h\in\mathcal{H}(\lambda,\alpha)$.
    Then there exists a homogeneous Moran set $C$ such that $\dimL C=\lambda$, $\dimA C = \alpha$, $\underline{\dim}_\theta C=h(0)$, and
    \begin{equation*}
        \overline{\dim}_\theta C=h(\theta)
    \end{equation*}
    for all $\theta\in[0,1]$.
\end{theorem}
\begin{proof}
    We will assume that $\lambda\leq h(0)<h(\theta)<\alpha$ for all $\theta\in(0,1]$, and that there exists $\theta_0 > 0$ such that $h(\theta_0) < h(1)$.
    The other cases are easier and can be proven with minor modifications or using a direct Moran construction.

    Let $(\epsilon_n)_{n=1}^\infty\subset(0,1)$ converge monotonically to 0 and $(\gamma_n)_{n=1}^\infty\subset(\lambda,\alpha)$ converge monotonically to $h(0)$ in such a way that $\gamma_n < h(\epsilon_{n+1}) < h(1)$ for all $n \in \N$, and $\gamma_n + 1/n \leq h(\epsilon_{n+1})$ for all $n$ sufficiently large.
    Note that if $h(0)>\lambda$ we can take $\gamma_n=h(0)$ for all $n$.
    We will define what we refer to as \emph{mountains} $f_n$ and \emph{valleys} $e_n$; a valley will be used to connect two mountains.
    Graphical representations of the functions $f_n$ and $e_n$ are given in \cref{f:mountain-ctr} and \cref{f:valley-ctr} respectively.
    Then we will define a function $g$ by concatenating the $f_n$ and $e_n$, and the corresponding Moran set $C$ will be given by \cref{l:Moran-formula}.

    The functions $f_n$ will ensure that $\overline{\dim}_\theta C=h(\theta)$ for $\theta>0$, and the functions $e_n$ will ensure that $g$ is continuous, $\dimH C=h(0)$, $\dimL C=\lambda$, and $\dimA C=\alpha$.

    \begin{proofpart}
        Construction of the mountain $f_n\colon[0,\log(1/\epsilon_n)]\to[h(\epsilon_n),h(1)]$ for $n\in\N$.
    \end{proofpart}
    \begin{figure}[t]
        \centering
        \def\xstar{0.994571}
\begin{tikzpicture}[>=stealth,scale=4,font=\tiny]
    \draw[->] (-0.1,0) -- (3.2,0);
    \draw[thick] plot file {figures/left_bump_points.txt};
    \draw[thick] plot file {figures/right_bump_points.txt};
    \draw[] (0.00,-0.05) node[below]{$0$} -- (0.00,0.05);
    \draw[] (3,-0.05) node[below]{$\log(1/\epsilon_n)$} -- (3,0.05);
    \draw[dashed] (0.319431,-0.05) node[below]{$x$} -- (0.319431,0.31875);
    \draw[dashed] (\xstar,-0.05) node[below]{$x^*$} -- (\xstar,0.6);
    \draw[dashed] (2.5,-0.05) node[below]{$x+\log(1/\theta)$} -- (2.5,0.31875);
    \draw[dotted]  (0,0.6) node[left]{$h(1)$}-- (3,0.6);
    \draw[dotted]  (0,0.31875) node[left]{$h(\theta)$}-- (3,0.31875);
    \draw[dotted]  (0,0.1) node[left]{$h(\epsilon_n)$}-- (3,0.1);
\end{tikzpicture}
        \caption{The construction of the \emph{mountain} $f_n$.}
        \label{f:mountain-ctr}
    \end{figure}
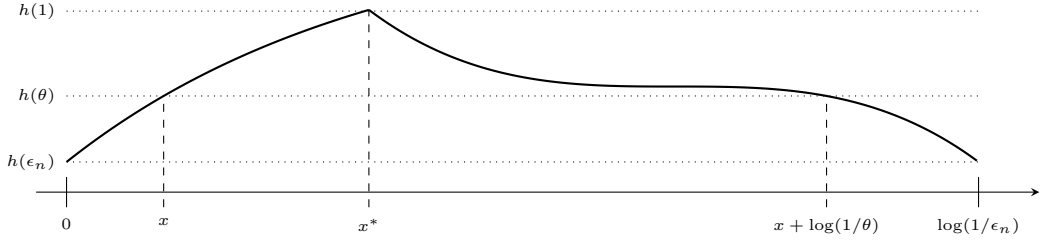
    First set
    \begin{equation*}
        x^*\coloneqq\log\left(\frac{\alpha-h(\epsilon_n)}{\alpha-h(1)}\right)
    \end{equation*}
    and for $x\in [0,x^*]$ define $f_n(x)=\alpha-(\alpha-h(\epsilon_n))\exp(-x)$.
    Observe that $f_n(0)=h(\epsilon_n)$, $f_n(x^*)=h(1)$, and
    \begin{equation}\label{e:0-star-bound}
        \dinil{-}f_n(x)=(\alpha-h(\epsilon_n))\exp(-x) = \alpha-f_n(x)
    \end{equation}
    for $x\in(0,x^*]$.

    Now for $x\in[0,x^*]$, if $\theta_x$ is such that $h(\theta_x)=f_n(x)$, we define $f_n(x+\log(1/\theta_x))=h(\theta_x)$.
    This is well-defined since $h$ is non-decreasing and continuous.
    In particular, $f_n(x)$ is non-increasing and continuous on $[x^*,\log(1/\epsilon_n)]$ with $f_n(\log(1/\epsilon_n))=f_n(0)=h(\epsilon_n)$.

    We now wish to bound $\dinil{-}f_n(x+\log(1/\theta_x))$ for $x\in(0,x^*]$.
    First, note that
    \begin{equation*}
        x=\log\left(\frac{\alpha-h(\epsilon_n)}{\alpha-h(\theta_x)}\right).
    \end{equation*}
    Then rearranging \cref{e:h-bound-general}, we obtain
    \begin{equation*}
        \diniu{+}h(\theta_x)\leq (h(\theta_x)-\lambda)\left(\frac{\diniu{+}h(\theta_x)}{h(\theta_x)-\alpha}+\frac{1}{\theta_x}\right).
    \end{equation*}
    Since $h(\theta_x)<\alpha$, $\diniu{+}h(\theta_x)<\frac{\alpha-h(\theta_x)}{\theta_x}$ so that $\frac{\diniu{+}h(\theta_x)}{h(\theta_x)-\alpha}+\frac{1}{\theta_x}>0$.
    Therefore,
    \begin{equation}\label{e:hbd-2}
        \frac{\diniu{+}h(\theta_x)}{\frac{\diniu{+}h(\theta_x)}{\alpha-h(\theta_x)}-\frac{1}{\theta_x}}\geq\lambda-h(\theta_x)=\lambda-f_n(x+\log(1/\theta_x)).
    \end{equation}
    But $x+\log(1/\theta_x)$ is a smooth function of $\theta_x$ and $h(\theta_x)$, and $h(\theta_x)=f_n(x+\log(1/\theta_x))$, so since $x$ decreases as $\theta_x$ decreases,
    \begin{equation*}
        \diniu{+}h(\theta_x)= \dinil{-}f_n(x+\log(1/\theta_x))\cdot\left(\frac{\diniu{+}h(\theta_x)}{\alpha-h(\theta_x)}-\frac{1}{\theta_x}\right)
    \end{equation*}
    which when combined with \cref{e:hbd-2} yields $\dinil{-}f_n(x+\log(1/\theta_x))\geq \lambda-f_n(x+\log(1/\theta_x))$.
    Note that we have shown that
    \begin{equation*}
        \dinil{-}f_n(x)\in[\lambda-f_n(x),\alpha-f_n(x)]
    \end{equation*}
    for all $x\in(0,\log(1/\epsilon_n)]$.

    \begin{proofpart}
        Construction of the valleys $e_n\colon[0,w_n]\to[\gamma_n,h(\epsilon_n)]$ where $w_n$ is given in \cref{e:wn-def} for $n\in\N$.
    \end{proofpart}
    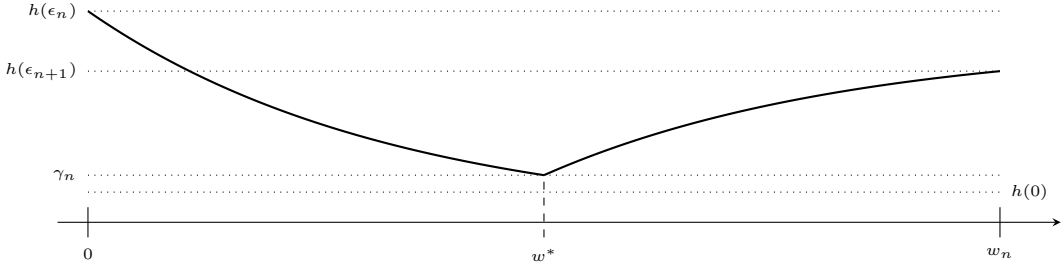
\begin{figure}[t]
        \centering
        \def\xstar{1.5}
\begin{tikzpicture}[>=stealth,scale=4,font=\tiny]
    \draw[->] (-0.1,0) -- (3.2,0);
    \draw[thick] plot file {figures/left_connector_points.txt};
    \draw[thick] plot file {figures/right_connector_points.txt};
    \draw[] (0,-0.05) node[below]{$0$} -- (0,0.05);
    \draw[] (3,-0.05) node[below]{$w_n$} -- (3,0.05);
    \draw[dashed] (\xstar,-0.05) node[below]{$w^*$} -- (\xstar,0.156191);
    \draw[dotted]  (0,0.7) node[left]{$h(\epsilon_n)$}-- (3,0.7);
    \draw[dotted]  (0,0.500973) node[left]{$h(\epsilon_{n+1})$}-- (3,0.500973);
    \draw[dotted]  (0,0.156171) node[left]{$\gamma_n$}-- (3,0.156191);
    \draw[dotted]  (0,0.1) -- (3,0.1)node[right]{$h(0)$};
\end{tikzpicture}
        \caption{The construction of the \emph{valley} $e_n$.}
        \label{f:valley-ctr}
    \end{figure}
    Set
    \begin{equation*}
        w^*\coloneqq\log\left(\frac{h(\epsilon_n)-\lambda}{\gamma_n-\lambda}\right)
    \end{equation*}
    and for $x\in[0,w^*]$ define $e_n(x)=\lambda-(\lambda-h(\epsilon_n))\exp(-x)$.
    Observe that $e_n(w^*)=\gamma_n$.
    Let
    \begin{equation}\label{e:wn-def}
        w_n\coloneqq w^*+\log\left(\frac{\alpha-\gamma_n}{\alpha-h(\epsilon_{n+1})}\right)
    \end{equation}
    and for $x\in[w^*,w_n]$ define $e_n(x)=\alpha-(\alpha-\gamma_n)\exp(-x+w^*)$.
    Of course, $e_n(w_n)=h(\epsilon_{n+1})$.
    It is clear that $\dinil{-}e_n(x)=\lambda-e_n(x)$ for $x\in(0,w^*]$ and $\dinil{-}e_n(x)=\alpha-e_n(x)$ for all $x\in(w^*,w_n]$.

    \begin{proofpart}
        Construction of $g\in\mathcal{G}(\lambda,\alpha)$ and the corresponding Moran set $C$.
    \end{proofpart}
    Let $\tilde{g}$ denote the concatenation of the sequence $(f_1,e_1,f_2,e_2,\ldots)$.
    By \cref{l:g-check}, $\tilde{g}$ satisfies the hypotheses of \cref{l:exact-construction} (note that $g$ is not rapidly decreasing since since $\limsup_{x\to\infty}\tilde{g}(x)>0$), and get a corresponding function $g$ and sequence $\bm{r}$.
    Note that $g\in\mathcal{G}(\lambda,\alpha)$.
    Let $C=C(\bm{r})$ denote the corresponding Moran set.

    That $\overline{\dim}_\theta C=h(\theta)$ for $\theta\in(0,1]$ follows by definition of the functions $f_n$ and the fact that
    \begin{equation*}
        \lim_{n\to\infty}\left(\sup_{x\in[0,w_n]}e_n(x)\right)\leq\lim_{\theta\to 0}h(\theta).
    \end{equation*}
    Moreover, \cref{l:Moran-formula} directly gives that $\dimH C=\underline{\dim}_\theta C = h(0)$ for $\theta\in[0,1]$, $\lambda\leq\dimL C$, and $\dimA C\leq\alpha$.

    To see that $\dimA C\geq\alpha$, note that the derivative of the strictly increasing part of each mountain is uniformly bounded above by $\alpha$. Therefore, for all $n \in \N$, the length of the domain of the $n$\textsuperscript{th} mountain can be uniformly bounded below:
    \begin{equation*}
        \log \left(\frac{1}{\epsilon_n}\right) \geq \log \left(\frac{1}{\epsilon_1}\right) > \frac{h(1)-h(\epsilon_1)}{\alpha} > 0.
    \end{equation*}
    Similarly, for all $n$ sufficiently large, the length $w_n$ of the domain of the $n$\textsuperscript{th} valley can be bounded below by $1/(\alpha n)$.
    Therefore there exists $\delta>0$ and a sequence $(b_m)_{m=1}^\infty$ such that for all $m \in \N$ we have $b_m \geq \delta m$, and $\diniu{+}g(x)+g(x)=\alpha$ for all $x\in[b_m,b_m+\delta/m]$.
    Then
    \begin{align*}
        \exp\left(b_m + \frac{\delta}{m}\right) - \exp(b_m) &=  \left(\exp\left(\frac{\delta}{m}\right) - 1\right) \cdot \exp(b_m)\\
                                                            &\geq \frac{\delta \cdot \exp(\delta m)}{m} \fto{m\to\infty} \infty.
    \end{align*}
    Thus we can define a function $\psi \colon \R \to \R^+$ such that $\psi(b_m) = b_m + \delta/m$ for all $m\in\N$ and $\lim_{x\to\infty}(\exp(\psi(x))-\exp(x))=\infty$.
    In particular, $D^+g(y) + g(y) = \alpha$ for all $y \in [b_m,\psi(b_m)]$ and infinitely many $m$.
    By \cref{l:Moran-formula}~\ref{im:A-dim-block}, it follows that $\dimA C\geq\alpha$.

    An analogous application of \cref{l:Moran-formula}~\ref{im:L-dim-block} gives that ${\dimL C\leq\lambda}$.
\end{proof}
\subsection{Prescribing the intermediate dimensions}\label{ss:pres}
We can get more varied behaviour for the lower intermediate dimensions by taking a finite union of Moran sets.
For example, the following proposition is straightforward to verify:
\begin{proposition}\label{p:finite-union}
    Suppose $g_i\in\mathcal{G}(0,d)$ for $i=1,\ldots,m$ have corresponding sequences $\bm{r}_i\subset(0,1/2]$ satisfying
    \begin{equation*}
        |g_i(x)-s_{\bm{r}_i}(\exp(-\exp(x)))|\leq d\log(2)\exp(-x).
    \end{equation*}
    Let $M$ be a disjoint union of translations of the homogeneous Moran sets $C(\bm{r}_i)$.
    Then for $\theta\in(0,1]$,
    \begin{enumerate}[nl]
        \item $\displaystyle\overline{\dim}_\theta M=\limsup_{x\to\infty}\max_{i=1,\ldots,m}\bigl(\inf_{y\in[x,x+\log(1/\theta)]}g_i(y)\bigr)$,
        \item $\displaystyle\underline{\dim}_\theta M=\liminf_{x\to\infty}\max_{i=1,\ldots,m}\bigl(\inf_{y\in[x,x+\log(1/\theta)]}g_i(y)\bigr)$,
        \item $\displaystyle\dimH M=\max_{i=1,\ldots,m}\liminf_{x\to\infty}g_i(x)$.
    \end{enumerate}
\end{proposition}
Suppose $h(\theta)$ satisfies $h(\epsilon)=h(0)$ for some $\epsilon>0$, and let $g$ denote the infinite concatenation of a mountain $f\colon[0,\log(1/\epsilon)]\to(0,d)$ constructed as in \cref{t:upper-h-form}.
If $C$ denotes the corresponding Moran set, then $\overline{\dim}_\theta C=h(\theta)$.

Now suppose $N$ is large, and define functions $g_i\coloneqq\kappa_{w_i}(g)$ where $w_i=\frac{(i-1)}{N}\log(1/\epsilon)$ for each $i\in\{1,\ldots,N\}$.
Write $A=d\log(1/\epsilon)$.
Then if $x$ is arbitrary, since the $g_i$ are Lipschitz continuous with constant $d$, there is some $i$ depending on $x$ such that
\begin{equation*}
    \inf_{y\in[x,x+\log(1/\theta)]}g_i(y)\geq h(\theta)-\frac{A}{N}
\end{equation*}
for all large $x$.
In particular, if $M$ denotes the set given by \cref{p:finite-union}, this implies that
\begin{equation*}
    h(\theta)-\frac{A}{N}\leq\underline{\dim}_\theta M\leq\overline{\dim}_\theta M= h(\theta).
\end{equation*}
In other words, by taking a finite union of homogeneous Moran sets, we can ensure that the upper and lower intermediate dimensions are arbitrarily close.

Motivated by this observation, we now construct a set such that the intermediate dimensions exist and are given by a prescribed formula $h(\theta)$.
At a fixed scale $\delta>0$, the set $M$ will look like a finite union of Moran sets each with the same upper intermediate dimensions.
As $\delta$ goes to zero, the resolution increases, so that the intermediate dimensions exist.
The construction here is mildly complicated by the fact that the mountains $f_n$ and valleys $e_n$ can have arbitrarily large support if $h(\theta)>h(0)$ for all $\theta>0$.
\begin{theorem}\label{t:gen-h-form}
    Let $0\leq\lambda\leq\alpha\leq d$ and let $h\in\mathcal{H}(\lambda,\alpha)$.
    Then there exists a compact perfect set $M$ such that $\dimL M=\lambda$, $\dimA M = \alpha$ and
    \begin{equation*}
        \dim_\theta C=h(\theta)
    \end{equation*}
    for all $\theta\in[0,1]$.
\end{theorem}
\begin{proof}
    We will assume that $\lambda\leq h(0)<h(\theta)<\alpha$ for all $\theta\in(0,1]$, and that there exists $\theta_0 > 0$ such that $h(\theta_0) < h(1)$. 
    The remaining cases follow by similar, but slightly easier, arguments.
    \begin{proofpart}
        Construction of the set $M$.
    \end{proofpart}
    As in the proof of \cref{t:upper-h-form}, fix non-increasing sequences $(\epsilon_n)_{n=1}^\infty$ and $(\gamma_n)_{n=1}^\infty$ and construct corresponding mountains $(f_n)_{n=1}^\infty$ defined on intervals $[0,z_n]$ and valleys $(g_n)_{n=1}^\infty$ defined on intervals $[0,w_n]$, where $z_n=\log(1/\epsilon_n)$ and $w_n$ is defined as in \cref{e:wn-def}.
    We may choose $\epsilon_n$ and $\gamma_n$ so that $w_n+z_n=2^n$.

    Let $\Psi=\{0,1\}\times\{0,1,2,3\}$ and let $\Psi^*=\bigcup_{n=0}^\infty\Psi^n$.
    We first associate to each $\eta\in\Psi^*$ a number $a(\eta)\in[0,\infty)$ as follows.
    Given $k\in\N$ and $i=(u,v)\in\Psi$, we define
    \begin{equation*}
        \psi(k,i)= u 2^{-k}+v 4^{k-1}
    \end{equation*}
    and then for $\eta=(i_1,\ldots,i_k)$, we set
    \begin{equation*}
        a(\eta)=\sum_{n=1}^k\psi(n,i_n).
    \end{equation*}
    Observe that $a(\Psi^k)=\{j2^{-k}:j\in\Z\}\cap[0,4^{k})$.

    For $k\in\N$ and $i\in\Psi$, we define $c_{k,i}(x)=h(\epsilon_k)$ for all $x\in[0,\psi(k,i)]$.
    Now for each $\eta=(i_1,\ldots,i_n)\in\Psi^*$, let $\tilde{g}_\eta$ denote the concatenation of the sequence
    \begin{equation*}
        (f_1,e_1,c_{1,i_1},f_2,e_2,c_{2,i_2},\ldots,f_n,e_n,c_{n,i_n},f_{n+1},e_{n+1},f_{n+2},e_{n+2},\ldots)
    \end{equation*}
    and set $g_\eta\coloneqq\kappa_{w_0}(\tilde{g}_\eta)$, where $w_0$ is guaranteed by \cref{l:exact-construction} and does not depend on the choice of $\eta$ since $g_\eta(0)=f_1(0)$ for all $\eta$, and moreover $w_0$ can be taken to be arbitrarily large.

    Thus there is a sequence $\bm{r}(\eta)\coloneqq(r_j(\eta))_{j=1}^\infty\subset(0,1/2]$ such that for all $x\geq w_0$,
    \begin{equation*}
        |s_\eta(\exp(-\exp(x)))-g_{\eta}(x)|\leq d\log(2)\cdot\exp(-x)
    \end{equation*}
    where $s_\eta\coloneqq s_{\bm{r}(\eta)}$.

    Let $\varnothing$ denote the word of length $0$, and let $\rho_k=r_1(\varnothing)\cdots r_k(\varnothing)$.
    For $k\geq 0$, let $y_k=w_0+\sum_{i=1}^k(w_i+z_i)=w_0+2^{k+1}-1$.
    Then let $n_k$ be the maximal index such that $\log\log(1/\rho_{n_k})\leq y_k$.
    Choosing $w_0$ large, we may assume that $n_k\geq 3k$ for all $k\in\N$.

    Let $\mathcal{I}=\{0,1\}^d$ and let $L\colon\mathcal{I}^3\to\Psi$ be given by $L(\bm{i},\bm{j},\bm{k})=(\bm{i}^{(1)},\bm{j}^{(1)}+2(\bm{k}^{(1)}))$.
    For $\ell\in\N$, we let $k_\ell$ denote the maximal index such that $n_{k_\ell}\leq\ell$.
    We then define a map $\Lambda\colon\mathcal{I}^*\to\Psi^*$ by
    \begin{equation*}
        \Lambda(\bm{i}_1,\ldots,\bm{i}_\ell) = (L(\bm{i}_1,\bm{i}_2,\bm{i}_3),L(\bm{i}_4,\bm{i}_5,\bm{i}_6),\ldots,L(\bm{i}_{3(k_\ell-1)+1},\bm{i}_{3(k_\ell-1)+2},\bm{i}_{3(k_\ell-1)+3})).
    \end{equation*}
    This is well-defined since $\ell\geq n_{k_\ell}\geq 3k_\ell$.

    We now construct our inhomogeneous Moran set $M$ as follows.
    Given a word $\sigma=(\bm{i}_1,\ldots,\bm{i}_\ell)\in\mathcal{I}^\ell$, let $\eta=\Lambda(\sigma)$.
    We then set $S_\sigma=S^1_{\bm{i}_1,\eta}\circ\cdots\circ S^\ell_{\bm{i}_\ell,\eta}$ where $S^i_{\bm{i},\eta}(x)=r_i(\eta)\cdot x+b^i_{\bm{i}}(\eta)$ with
    \begin{equation*}
        b^i_{\bm{i}}(\eta)^{(j)} =
        \begin{cases}
            0 &: \bm{i}^{(j)}=0\\*
            1-r_i(\eta) &: \bm{i}^{(j)}=1
        \end{cases}
        .
    \end{equation*}
    We now set
    \begin{equation*}
        M_\ell=\bigcup_{\sigma\in\mathcal{I}^\ell}S_\sigma([0,1]^d).
    \end{equation*}
    Note that if $\sigma$ is a prefix of $\tau$, then $\Lambda(\sigma)$ is a prefix of $\Lambda(\tau)$ and therefore $S_\sigma([0,1]^d)\supseteq S_\tau([0,1]^d)$.
    Thus $M_0\supseteq M_1\supseteq\cdots$ so that the set
    \begin{equation*}
        M\coloneqq\bigcap_{\ell=0}^\infty M_\ell
    \end{equation*}
    is non-empty.

    Intuitively, at a fixed scale $\delta$, $M$ looks like a union of $8^k$ homogeneous Moran sets corresponding to the sequences $\bm{r}(\eta)$ for $\eta\in\Psi^k$.
    We can make this precise in the following sense.
    For $\eta\in\Psi^k$, we define
    \begin{align*}
        \mathcal{B}_k(\eta) &= \{(\sigma_1,\ldots,\sigma_k)\in\mathcal{I}^{3k}:L(\sigma_i)=\eta_i\text{ for each }1\leq i\leq k\},\\
        J_\eta &= \bigcup_{\sigma\in\mathcal{B}_k(\eta)}S_\sigma([0,1]^d).
    \end{align*}
    Let $C_\eta\coloneqq C(\bm{r}(\eta))$ denote the homogeneous Moran set corresponding to the function $g_\eta$.
    Let $\ell\in\N$ satisfy
    \begin{equation*}
        y_k+a(\eta^-)<\log\log(1/(r_1(\eta)\cdots r_\ell(\eta)))\leq y_{k+1}+a(\eta)
    \end{equation*}
    where $\eta^-\in\Psi^{k-1}$ is the unique prefix of $\eta$.
    Since $g_\varnothing(y_{k})=g_\eta(y_{k}+a(\eta^-))$, if $\sigma\in \mathcal{I}^\ell$, then $\eta$ is a prefix of $\Lambda(\sigma)$.
    Moreover, if $\tau\in\Psi^*$ is any word with $\eta$ as a prefix, $g_\tau(x)=g_\eta(x)$ for all $x\leq y_{k+1}+a(\eta)$.
    Thus for any such $\ell$, we have
    \begin{equation}\label{e:Ml-formula}
        M_\ell\cap J_\eta=(C_\eta)_\ell\cap J_\eta.
    \end{equation}
    But then if $\eta'$ is a prefix of $\eta$, $r_\ell(\eta')=r_\ell(\eta)$ for all $\ell$ such that
    \begin{equation}\label{e:eta-valid}
        \log\log(1/(r_1(\eta)\cdots r_\ell(\eta)))\leq y_{k+1}+a(\eta).
    \end{equation}
    Thus \cref{e:Ml-formula} holds for any $\ell$ satisfying \cref{e:eta-valid}.

    We also note that $(C_\eta)_\ell\cap J_\eta$ consists of exactly $2^{d\ell-3k}$ hypercubes with diameter $r_1(\eta)\cdots r_\ell(\eta)$.
    \begin{proofpart}
        Proof that $\dim_\theta M=h(\theta)$ for $\theta\in(0,1]$.
    \end{proofpart}

    Fix $\theta\in(0,1]$.
    We first show that $\overline{\dim}_\theta M\leq h(\theta)$.
    Let $\delta$ be sufficiently small so that $\delta\leq\rho_{k_0}$ where $\epsilon_{k_0}\leq\theta$.
    Now let $k$ be such that $\rho_{n_{k}}<\delta^{1/\theta}$.
    It now follows by the same argument as \cref{l:flat-covers} that for each $\eta\in\Psi^k$, with $s_\eta\coloneqq\inf_{\phi\in[\delta^{1/\theta},\delta]}s_\eta(\phi)$,
    \begin{align*}
        \inf\Bigl\{\mathcal{C}^{s_\eta}(\mathcal{U}):\mathcal{U}\text{ is a $(\delta,\theta)$-cover of $(C_\eta)_{\ell(\eta)}\cap J_\eta$}\Bigr\}\leq 8^{-k}
    \end{align*}
    where $\ell(\eta)$ is minimal such that $r_1(\eta)\cdots r_\ell(\eta)\leq\delta^{1/\theta}$.
    But $\ell(\eta)$ satisfies \cref{e:eta-valid} since $\rho_{n_{k}}<\delta^{1/\theta}$, so that $M\subseteq \bigcup_{\eta\in\Psi^k}(C_\eta)_\ell\cap J_\eta$.
    Therefore, $s_\eta\leq h(\theta)+d\log(2)\cdot\exp(-y_{n_k})$.
    This implies that $\overline{\dim}_\theta M\leq h(\theta)$.

    Now fix $\epsilon>0$: we will show that $\underline{\dim}_\theta M\geq h(\theta)-(2+d)\epsilon$.
    The various variables in this proof are depicted in \cref{f:offset-diag}.
    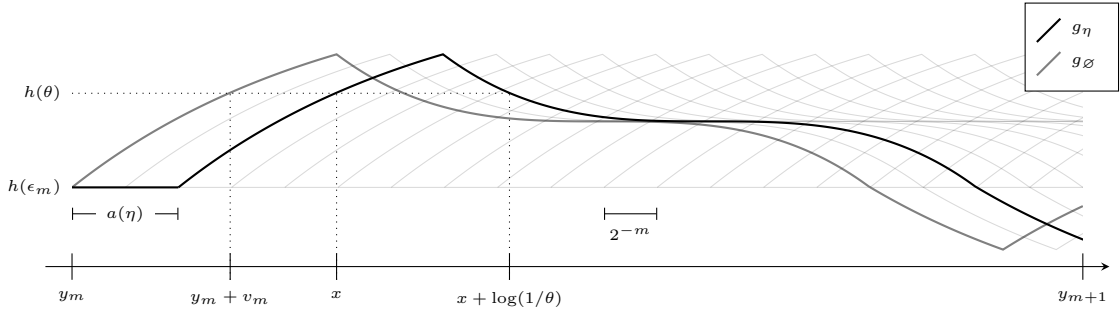
\begin{figure}[t]
        \centering
        \def\ymvm{0.594571}
\def\xval{0.994571}
\def\xtval{1.64457}
\def\htheta{0.656538}
\begin{tikzpicture}[>=stealth,scale=3.5,font=\tiny]
    \draw[->] (-0.1,0) -- (3.9,0);
    \begin{scope}
        \clip (0,0) rectangle (3.8,1);
        \draw[thin,opacity=0.15] (0,0.3) -- (3.8, 0.3);

        \draw[thick,gray] plot file {figures/offset_points.txt}; %
        \draw[thick,xshift=0.4cm] plot file {figures/offset_points.txt}; %
        \draw[thick] (0,0.3) -- (0.4, 0.3);

        \foreach\shift in {0.2,0.4,...,4} {
            \draw[thin,opacity=0.15,xshift=\shift cm] plot file {figures/offset_points.txt};
        }
    \end{scope}
    \draw[] (0.00,-0.05) node[below]{$y_m$} -- (0.00,0.05);
    \draw[] (3.8,-0.05) node[below]{$y_{m+1}$} -- (3.8,0.05);
    \draw[] (\ymvm,-0.05) node[below]{$y_m+v_m$} -- (\ymvm,0.05);
    \draw[] (\xval,-0.05) node[below]{$x$} -- (\xval,0.05);
    \draw[] (\xtval,-0.05) node[below]{$x+\log(1/\theta)$} -- (\xtval,0.05);

    \draw[thin,dotted] (\ymvm,-0.05) -- (\ymvm,\htheta);
    \draw[thin,dotted] (\xval,-0.05) -- (\xval,\htheta);
    \draw[thin,dotted] (\xtval,-0.05) -- (\xtval,\htheta);

    \draw[thin,dotted] (0,\htheta) node[left]{$h(\theta)$}-- (\xtval,\htheta);

    \draw[dotted] (0,0.3) node[left]{$h(\epsilon_m)$}-- (0,0.3);
    \draw[|-|] (0,0.2) -- node[fill=white]{$a(\eta)$} (0.4,0.2);
    \draw[|-|] (2,0.2) -- node[below]{$2^{-m}$} (2.2,0.2);

    \matrix [draw, fill=white, below left] at (current bounding box.north east) {
        \node [strike out,draw=black,thick,label=right:{$g_\eta$}] {}; \\
        \node [strike out,draw=gray,thick,label=right:{$g_\varnothing$}] {}; \\
    };
\end{tikzpicture}
        \caption{Choice of $g_\eta$ for the lower bound of $\dim_\theta M$.}
        \label{f:offset-diag}
    \end{figure}
    Let $k$ be such that $2^{-k}\leq\epsilon$.
    Let $\delta>0$ be small and let $x\coloneqq\log\log(1/\delta)$.
    We may assume that
    \begin{enumerate}[nl,a]
        \item $d\log(2)\exp(-x)\leq\epsilon$,
        \item $x\geq y_k$, and
        \item $x\geq y_m$ for some $m$ with $\epsilon_m\leq\theta$.
    \end{enumerate}

    For each $m\in\N$, there is some $v_m$ such that $f_m(v_m)=h(\theta)$.
    Equivalently, $g_\varnothing(y_m+v_m)=h(\theta)$.
    Let $m$ be maximal such that $y_m+v_m\leq x$.
    Since $y_{m+1}+v_{m+1}-(y_m+v_m)\leq 4^m$, there is some $\eta_0\in\Psi^m$ such that $|a(\eta_0)-(x-y_m-v_m)|\leq 2^{-k}$.
    Then since $\diniu{+}g_{\eta}(x)\in[-d,d]$ for all $x\in\R$,
    \begin{equation*}
        \inf_{\phi\in[\delta^{1/\theta},\delta]}s_\eta(\phi)\geq \inf_{y\in[x,x+\log(1/\theta)]}g_\eta(y)-\epsilon\geq h(\theta)-(1+d)\epsilon.
    \end{equation*}
    Set $s=h(\theta)-(2+d)\epsilon$.
    Again by the same argument as \cref{l:flat-covers}, since $x+\log(1/\theta)<y_{m+1}+v_{m+1}< y_{m+2}$, with $\eta\in\Psi^{m+1}$ satisfying $g_{\eta_0}=g_\eta$, we have
    \begin{align*}
        C\cdot\frac{\delta^{-\epsilon}}{8^m}&\leq \inf\Bigl\{\mathcal{C}^s(\mathcal{U}):\mathcal{U}\text{ is a $(\delta,\theta)$-cover of $M\cap J_\eta$}\Bigr\}\\
                                                        &\leq \inf\Bigl\{\mathcal{C}^s(\mathcal{U}):\mathcal{U}\text{ is a $(\delta,\theta)$-cover of $M$}\Bigr\}
    \end{align*}
    for some constant $C>0$ not depending on $\delta$.
    But $x\geq y_m\geq 2^m-1$, so
    \begin{equation*}
        \frac{\delta^{-\epsilon}}{8^m}\geq\frac{\bigl(\exp(\exp(2^m-1))\bigr)^\epsilon}{8^m}\fto{m\to\infty}\infty
    \end{equation*}
    as required.
    \begin{proofpart}
        Proof that $\dimH M=h(0)$, $\dimL M=\lambda$, and $\dimA M=\alpha$.
    \end{proofpart}
    It is clear that $\dimH M\geq h(0)$ since $\liminf_{\delta\to 0}s_\eta(\delta)\geq h(0)$ for all $\eta\in\Psi^*$.
    Conversely, let $\epsilon>0$: we will show that $\dimH M\leq h(0)+2\epsilon$.
    Let $n_0$ be sufficiently large so that $\gamma_{n_0}\leq h(0)+\epsilon$.
    Then let $\delta>0$ be sufficiently small so that with $x=\log\log(1/\delta)$, we have $x\geq y_{n_0+1}$ and $d\log(2)\cdot\exp(-x)\leq\epsilon$.

    Let $m$ be such that $x<y_m$.
    For each $\eta\in\Psi^m$, by the choice of $n_0$, there exists some $x\leq x_\eta<y_{m+1}+a(\eta)$ such that
    \begin{equation*}
        g_\eta(x_\eta)= \gamma_m\leq h(0)+\epsilon.
    \end{equation*}
    Then by the same argument as \cref{l:flat-covers}, since $d\log(2)\cdot\exp(-x)\leq\epsilon$, with $s=h(0)+2\epsilon$ we have with $\ell(\eta)$ minimal so that $\log\log(1/\rho_{\ell(\eta)})\geq x_\eta$,
    \begin{align*}
        \inf\Bigl\{\mathcal{C}^s(\mathcal{U}):\mathcal{U}\text{ is a $(\delta,0)$-cover of $(C_\eta)_{\ell(\eta)}\cap J_\eta$}\Bigr\}\leq 8^{-m}.
    \end{align*}
    Moreover, for $x$ sufficiently large, we can ensure $\log\log(1/\rho_{\ell(\eta)})\leq y_{m+1}\leq y_{m+1}+a(\eta)$.
    Thus $M\subseteq \bigcup_{\eta\in\Psi^m}(C_\eta)_{\ell(\eta)}\cap J_\eta$ so that
    \begin{equation*}
        \inf\Bigl\{\mathcal{C}^s(\mathcal{U}):\mathcal{U}\text{ is a $(\delta,0)$-cover of $M$}\Bigr\}\leq 1.
    \end{equation*}
    But $\delta>0$ was arbitrary, so that $\dimH M\leq h(0)+2\epsilon$, as required.

    Now we will see that $\dimA M=\alpha$; the proof that $\dimL M=\lambda$ follows similarly.
    To see that $\dimA M\geq\alpha$, observe that there is some $\delta>0$ such that $\diniu{+}g_\varnothing(x)+g_\varnothing(x)=\alpha$ for all $m\in\N$ and $x\in[y_{m+1}-\delta,y_{m+1}]$.
    Let $\tau=\{(0,0),(0,0,),\ldots,(0,0)\}\in\Psi^m$ and observe that $g_\varnothing=g_\tau$.
    Then if $\log\log(1/\rho_\ell)\in[y_{m+1}-\delta,y_{m+1}]$, we have $M_\ell\cap J_\tau = C_\varnothing\cap J_\tau$.
    Thus $\dimA M\geq\alpha$ follows by the same computation from \cref{t:upper-h-form}.

    Conversely, it suffices to show for all $\epsilon>0$ there is some $a>0$ so that for all $\ell\in\N$, $I\in M_\ell$ with $\diam I=R$,
    \begin{equation*}
        N_r(I\cap M)\leq (R/r)^{\alpha+2\epsilon}
    \end{equation*}
    for all $0<r\leq aR$.
    Let $m,k$ be minimal so that $\log\log(1/r)\leq\log\log(1/\rho_m)\leq y_k$.

    First suppose $\ell\geq 3k$.
    Then there is a unique $\eta\in\Psi^k$ such that $I\subset J_\eta$, so that
    \begin{equation*}
        (C_\eta)_j\cap I\cap M_j=I\cap M_j
    \end{equation*}
    for all $\ell\leq j\leq m$.
    Then since $\diniu{+}g_\eta+g_\eta\leq\alpha$, a similar computation to the proof of \cref{l:Moran-formula} gives that for all $\epsilon>0$ and $\ell$ sufficiently large depending on $\epsilon$,
    \begin{equation*}
        N_r(I\cap M)\leq (2^d)^{m-\ell}\leq\left(r_{\ell+1}(\eta)\cdots r_{m}(\eta)\right)^{-\bigl(\alpha+\frac{2d}{m-\ell}+\epsilon\bigr)}\leq(R/r)^{\alpha+\frac{2d}{m-\ell}+\epsilon}.
    \end{equation*}

    Otherwise, $\ell<3k$.
    Let $\eta\in\Psi^k$ satisfy $J_\eta\cap I\neq\varnothing$ and let $\sigma\in\mathcal{I}^{3k}$ have $S_\sigma([0,1]^d)\subseteq I\cap J_\eta$.
    Again,
    \begin{equation*}
        N_r(I\cap(C_\eta)_m)\leq(2^d)^{m-\ell}\leq(R/r)^{\alpha+\frac{2d}{m-\ell}}
    \end{equation*}
    so that
    \begin{align*}
        N_r(I\cap(C_\eta)_m\cap S_\sigma([0,1])^d)&\leq (2^d)^{m-3k}=(2^d)^{\ell-3k}(2^d)^{m-\ell}\\
                                               &\leq (2^d)^{\ell-3k}(R/r)^{\alpha+\frac{2d}{m-\ell}}.
    \end{align*}
    But $I\cap(C_\eta)_m\cap S_\sigma([0,1]^d)=I\cap M_m\cap S_\sigma([0,1]^d)$ and there are precisely $(2^d)^{3k-\ell}$ words $\sigma$, so that
    \begin{equation*}
        N_r(I\cap M)\leq N_r(I\cap M_m)\leq (R/r)^{\alpha+\frac{2d}{m-\ell}}.
    \end{equation*}
    We can therefore choose $a$ small enough so that in either case $N_r(I\cap M)\leq (R/r)^{\alpha+2 \epsilon}$, as required.
\end{proof}
\begin{remark}
    It is clear from the construction that there are inhomogeneous Moran sets for which any cover approximating the intermediate dimensions arbitrarily closely would require an unbounded number of scales as $\delta$ tends to zero.
    This answers a question of Falconer \cite{fal2021}.
\end{remark}
Using this construction, along with the preceding construction for the upper intermediate dimensions, we can now simultaneously prescribe the upper and lower intermediate dimensions.
\begin{corollary}\label{c:upper-lower-match}
    Let $0\leq\lambda\leq\alpha\leq d$ and let $\underline{h},\overline{h}\in\mathcal{H}(\lambda,\alpha)$ satisfy $\underline{h}(0)=\overline{h}(0)$ and $\underline{h}\leq\overline{h}$.
    Then there exists a compact perfect set $M$ such that $\dimL M=\lambda$, $\dimA M = \alpha$ and
    \begin{align*}
        \underline{\dim}_\theta C&=\underline{h}(\theta) & \overline{\dim}_\theta C&=\overline{h}(\theta)
    \end{align*}
    for all $\theta\in[0,1]$.
\end{corollary}
\begin{proof}
    Let $E,F$ be disjoint compact perfect sets such that $\dimL E=\dimL F=\lambda$, $\dimA E=\dimA F=\alpha$, $\dimH E=\dimH F=\underline{h}(0)=\overline{h}(0)$, and for $\theta\in(0,1]$
    \begin{equation*}
        \underline{\dim}_\theta F\leq\dim_\theta E=\underline{h}(\theta)\leq\overline{h}(\theta)=\overline{\dim}_\theta F.
    \end{equation*}
    For example, such a set $E$ is provided by \cref{t:gen-h-form} and such a set $F$ is provided by \cref{t:upper-h-form}.
    Let $M=E\cup F$.

    Then $\dimL M=\min\{\dimL E,\dimL F\}=\lambda$, $\dimA M=\max\{\dimA E,\dimA F\}=\alpha$,
    \begin{equation*}
        \underline{h}(\theta)=\underline{\dim}_\theta E\leq \underline{\dim}_\theta M\leq\max\{\overline{\dim}_\theta E,\underline{\dim}_\theta F\}=\underline{h}(\theta),
    \end{equation*}
    and
    \begin{equation*}
        \overline{\dim}_\theta M=\max\{\overline{\dim}_\theta E,\overline{\dim}_\theta F\}=\overline{h}(\theta)
    \end{equation*}
    for $\theta\in(0,1]$.
    Thus $M$ satisfies the requirements.
\end{proof}
\begin{acknowledgements}
    We thank Jonathan Fraser, Kenneth Falconer, Haipeng Chen, Lars Olsen and P\'eter Varj\'u for helpful comments.
    We also thank Kathryn Hare for pointing out some references on Assouad dimensions, and Boyuan Zhao for pointing out the reference for Dini derivatives.
    Finally, we thank the anonymous referee for several useful comments.

    AB was supported by a Leverhulme Trust Research Project Grant (RPG-2019-034).
    AR was supported by EPSRC Grant EP/V520123/1.
\end{acknowledgements}
\end{document}